\documentclass[11pt]{article}
\usepackage[utf8]{inputenc}
\usepackage{xcolor}
\usepackage{graphicx}
\usepackage{float}

\usepackage{amsmath,amsfonts,amssymb,mathrsfs}
\usepackage{amsthm}
\usepackage{graphicx}
\usepackage{tikz}
\usepackage{url}
\usepackage[colorlinks=true, linkcolor=blue, citecolor=blue, urlcolor=blue]{hyperref}
\usepackage{titling}
\usepackage{subcaption}

\usepackage{algorithm}
\usepackage{algorithmic}       

\usepackage{natbib}
\setcitestyle{citesep={;\allowbreak\space}}

\usepackage{titlesec}

\titlespacing*{\section}  {0pt}{.6ex plus .2ex minus .2ex}{.6ex}
\titlespacing*{\subsection}{0pt}{.4ex plus .15ex}{.4ex}

\usepackage{enumitem}

\setlist[itemize]{
  topsep   = 1pt,  
  partopsep= 0pt,  
  itemsep  = 1pt,  
  parsep   = 0pt   
}

\usepackage{booktabs}
\usepackage{tabularx}        
\usepackage{longtable}       
\usepackage{threeparttablex} 
\usepackage{siunitx}         

\newcolumntype{Y}{>{\raggedright\arraybackslash}X} 

\usepackage[english]{babel}
\newtheorem{theorem}{Theorem}[section]
\newtheorem{corollary}[theorem]{Corollary}
\newtheorem{lemma}[theorem]{Lemma}
\newtheorem{remark}{Remark}

\newtheorem{definition}{Definition}[section]
\newtheorem{assumption}{Assumption}[section]

\theoremstyle{definition}
\newtheorem{example}{Example}[section]

\numberwithin{equation}{section}

\DeclareMathOperator*{\esssup}{ess\,sup}
\DeclareMathOperator*{\essinf}{ess\,inf}

\usepackage{caption}
\usepackage{subcaption}

\usepackage{booktabs}        
\usepackage{tabularx}        
\usepackage{threeparttablex} 

\usepackage{tikz}
\usetikzlibrary{positioning}

\newcommand{\Tr}{\mathrm{Tr}}

\ifpdf
\hypersetup{
  pdftitle={Agency Problems and Adversarial Bilevel Optimization under Uncertainty and Cyber Threats},
  pdfauthor={Thibaut Mastrolia, Haoze Yan}
}
\fi

\title{Agency Problems and Adversarial Bilevel Optimization under Uncertainty and Cyber Threats}
\author{Thibaut \textsc{MASTROLIA} and Haoze \textsc{YAN}\thanks{\texttt{mastrolia@berkeley.edu}, \texttt{haoze.yan@berkeley.edu}. The authors thank 
Tom Bielecki, Caroline Hillairet, Emma Hubert, Dylan Possama\"i, Wissal Sabbagh for insightful discussion and feedback on preliminary results of this work. This work benefits from the France-Berkeley Fund 2023 as part of the project "Cyber risk mathematical modeling".}\\ UC Berkeley\\
Department of Industrial Engineering and Operations Research}

\begin{document}
\setlength{\droptitle}{-5cm}
\maketitle

\begin{abstract}
We study an agency problem between a leader (the principal) seeking to design an optimal incentive scheme to a follower (the agent) to increase the value of a risky project subjected to accidents and volatility uncertainty. The agency problem is formulated as a max-min bilevel stochastic control problem with accidents and ambiguity. We show that the problem of the follower is reduced to solve a second order BSDE with jumps, reducing the problem of the leader to solve an integro-partial Hamilton–Jacobi–Bellman–Isaacs (HJBI) equation.
{ By extending stochastic Perron's method to our setting, we obtain viscosity
sub- and supersolution envelopes for the Principal's integro-HJBI equation.
Under an additional comparison principle in a suitable polynomial-growth class,
these envelopes coincide and the Principal's value is identified with the
unique viscosity solution.}
The holding company seeks to design an optimal incentive scheme to mitigate these losses. In response, the subsidiary selects an optimal cybersecurity investment strategy, modeled through a stochastic epidemiological SIR (Susceptible-Infected-Recovered) framework. The cyber threat landscape is captured through an L-hop risk framework with two primary sources of risk, internal risk propagation via the contagion parameters in the SIR model, and external cyberattacks from a malicious external hacker. The uncertainty and adversarial nature of the hacking lead to consider a robust stochastic control approach that allows for increased volatility and ambiguity induced by cyber incidents. We illustrate our results with numerical simulations showing how the contracting mechanism enhances the quality of a cluster under cyber threats.
\end{abstract}

\section{Introduction}
According to Governor Michael S. Barr, speaking at the Federal Reserve Bank of New York on April 17, 2025 \textit{``Cybercrime is on the rise, and cybercriminals are increasingly turning to Gen AI to facilitate their crimes. Criminal tactics are becoming more sophisticated and available to a broader range of criminals. Estimates of direct and indirect costs of cyber incidents range from 1 to 10 percent of global GDP. Deepfake attacks have seen a twenty-fold increase in the last three years''.} Governor Barr’s remarks underscore the growing severity of cyber threats fueled by the hyper-connectivity of modern society. Individuals, businesses, public institutions, and critical infrastructure are increasingly interconnected through digital networks—creating vulnerabilities across virtually every sector. From social media platforms and private messaging services to healthcare systems, governments, and financial institutions, no domain is immune. These threats are not geographically confined either; cyberattacks are now a global concern, affecting nations and industries worldwide. Recent geopolitical developments—such as the Russia-Ukraine war—have further intensified cyber threats, particularly across Europe and NATO member states. Likewise, the COVID-19 pandemic, which accelerated the digitalization of services and online interaction, has expanded the attack surface for cybercriminals. However, cyber threats have been growing increasingly sophisticated over the past few decades, making it urgent to develop a strong agenda to address it as one of the main challenge of the 21st century (see, e.g., \cite{tatar2014comparative,karabacak2014strategies,eling2021cyber,amin2019practical,ghadge2020managing}). To address these challenges, the U.S. Department of Homeland Security’s Science and Technology Directorate has launched the Cyber Risk Economics (CyRiE) project. This initiative promotes research into the legal, behavioral, technical, and economic dimensions of cybersecurity. A key component of CyRiE focuses on designing effective incentives to optimize cyber-risk management, aiming to guide organizations in allocating resources toward the most impactful and valuable defenses. 

This work contributes to that objective by exploring how a parent (holding) firm can design optimal incentives and compensation mechanisms for its subsidiaries operating under cyber threat conditions. The goal is to ensure efficient monitoring and management of both the subsidiary’s portfolio and its cybersecurity strategies. 

\subsection{Incentives and agency theory}
Turning now to incentive mechanism, it has been investigated since the 1960s in economy and known as contract theory or agency problem, model with a Principal-Agent framework with information asymmetry. Holmstrom and Milgrom's 1987 pioneer work \cite{holmstrom1987aggregation} has set the paradigm in a continuous-time framework with continuous controlled process. It has then regained interest in the mathematical community in the last decades with the work of Sannikov \cite{sannikov2008continuous} and Cvitanic, Possamai and Touzi \cite{cvitanic2018dynamic,cvitanic2017moral}. In our model, the holding firm (the principal) monitors indirectly the action of the subsidiary (the agent) by proposing a compensation for its activities. The holding firm does not have a direct access to the activities of its subsidiary and only observes the result of its work through its wealth and corrupted devices in the SIR system. This asymmetry of information arises in a moral hazard situation in which the principal must anticipates the best reaction of the agent to propose an optimal incentives scheme. This problem is equivalent to solve a Stackelberg game in continuous time, see for example \cite{li2017review,hernandez2024principal,hernandez2024closed}. We usually address this problem as a bilevel stochastic optimization, in which the problem of the agent is embedded into the problem of the principal, known as \textit{the incentive compatibility condition} of the compensation offers by the principal to the agent ensuring the existence of a best reaction activity, see e.g. \cite{mastrolia2025agency,dempe2020bilevel}. We refer to \cite{tirole2010theory,cvitanic2012contract} for a more detailed overview of principal-agent, Stackelberg games and agency problem.

\paragraph{Stochastic control contributions.} The bilevel optimization investigated is
$$
{\
\begin{aligned}
V_0^P
:=&\sup_{\xi\in\Xi}\sup_{\hat\alpha(\xi)} \inf_{(\mathbb P,\eta)\in \mathcal H^{(\widehat\alpha(\xi))}(0,x_0)}
\mathbb E^{\mathbb P}\Big[F^P(X_T)-\xi-\int_0^T C^P\big(s,X_s,\widehat\alpha_s(\xi),\eta_s\big)ds
\Big]
\\
\text{subject to }\quad
&\text{(IC)}\ \ \widehat\alpha(\xi)\ \text{is an agent best response in \eqref{eq:agent-value}},\\
&\text{(IR)}\ \ V_0^A(\xi)\ \ge R_0,
\end{aligned}}
$$
where the risky project is solution to the following equation
$$
dX_t=b(t,X_{t-},\alpha_t,\eta_t)dt
   +\sigma(t,X_{t-},\eta_t)dW^{(\alpha,\eta)}_t
   +\int_E \beta(t,X_{t-},e)\tilde\mu^{(\alpha,\eta)}(dt,de),
$$
driven by the agent's control $\alpha$ under volatility uncertainty $\eta$ leading to an uncertain family of controlled probability $\mathcal P^{\alpha}$. 
Our work is the first one proposing (i) an applications to second order BSDE with jumps to stochastic control and volatility ambiguity resolving the agent's problem \eqref{eq:agent-value} below embedded in the leader-follower problem; (ii) extending stochastic Perron's method to stochastic control and max-min optimization with volatility uncertainty and jumps; (iii) developing a self-contained framework tractable for diverse applications including cyber risk management.

\paragraph{Cyber risk and L-hop propagation under ambiguity}
Cyberattacks vary widely in form and mechanism (see, e.g., \cite{uma2013survey,hathaway2012law,hillairet2023expansion,grove2019cybersecurity,boumezoued2023cyber,hillairet2024optimal}), but L-hop propagation models are particularly useful for capturing the dynamics of both external and internal threats. The term \textbf{L-hop} refers to the number of network connections (or "hops") an attack can traverse before reaching its target. External threats originate outside the network such as direct hacking attempts, modeled using a point process with exogenous intensity. Internal threats emerge from within the network, typically through infected nodes spreading malware or viruses. These internal dynamics are modeled using compartmental epidemiological models, such as the SIR (Susceptible-Infectious-Recovered) framework, see e.g. \cite{capasso1993mathematical,britton2010stochastic,elie2020contact}, in the context of cyber risk (see, e.g., \cite{del2022computational,hillairet2022cyber,hillairet2024optimal}). By integrating these components, the proposed model offers a robust framework for evaluating how financial firms can design efficient intra-organizational incentives that align cybersecurity investments with the broader objectives of risk mitigation and financial resilience. 

In the realm of cybersecurity, the inherent unpredictability and knowledge gaps that arise when constructing and deploying models to predict or prevent cyber-threats lead to various types of uncertainty. These uncertainties can arise from multiple sources and understanding them is vital for the development of more resilient and adaptive cybersecurity systems. This work focuses on three key types of uncertainty: (1) the propagation of cyber risk within the subsidiary cluster; (2) the impact on the system's wealth; and (3) the randomness and ambiguity inherent in the behavior of cyber attacks. This section introduces informally the problem investigated. A more rigorous framework is provided hereafter. 

As discussed previously, the propagation of a cyber attack is modeled using an epidemiological framework with stochastic noise. Specifically, we assume that the spread of the attack within the cluster, referred to as the internal L-hop risk, is governed by the following SIR (Susceptible-Infected-Recovered) system:
\begin{equation*}
    \begin{cases}
        dS_t = (-\bar\beta S_tI_t - \alpha_tS_t-\eta_tS_t)dt - \tilde\sigma(t,\eta_t) S_tI_td\widetilde W_t\\
        dI_t = (\bar\beta S_tI_t-\rho I_t +\eta_tS_t)dt + \tilde\sigma(t,\eta_t) S_tI_td\widetilde W_t\\
        dR_t = \rho I_tdt + \alpha_tS_tdt,
    \end{cases}
\end{equation*}
where $\eta$ denotes the unknown cyber attack and $\alpha$ the protection strategy used by the subsidiary. 

{Note that the classical deterministic SIR model does not contain a Brownian perturbation. 
In our stochastic cyber-risk model, the parameter $\bar\beta>0$ remains the baseline 
transmission rate used in the formal model of Section 5. The additional Brownian term 
models uncertainty in the cumulative infection pressure generated by cyber propagation. 
Equivalently, one may write the infinitesimal contagion pressure as
$$
dB_t=\bar\beta\,dt+\tilde\sigma(t,\eta_t)\,d\widetilde W_t,
$$
so that the stochastic contagion flow is $S_t I_t\,dB_t$. See \cite{hansen2011optimal,gray2011stochastic,hubert2022incentives}. Hence the susceptible and 
infected equations contain the stochastic terms $
-S_t I_t\,dB_t, +S_t I_t\,dB_t,$
respectively, together with the protection, recovered, and adversarial infection terms. 
Thus the randomness is not introduced through a separate transmission parameter; it is 
the diffusion component of the same SIR propagation mechanism, with volatility modulated 
by the cyber-attack variable $\eta$.} {
In our model, $\eta$ represents ambient, background cyber risk and volatility uncertainty, such as continuous vulnerability scanning, phishing campaigns, or the gradual degradation of IT infrastructure by stealthy actors. Because this type of threat operates in the background and its exact intensity is ambiguous, modeling it as a continuous factor affecting the diffusion/volatility is highly appropriate and aligns with continuous-time robust control frameworks. 
}

Regarding the uncertainty in the wealth of the subsidiary, we assume that the portfolio of the firm is given at time $t$ by the solution to the following SDE \begin{equation*}
     dP_s = P_s\left(\mu(s,I_s)dt + \sigma(s,I_s,\eta_s)dW_s + \int_{E}l_s(e)\mu_P(de,ds)\right),
\end{equation*}
where $\mu_P$ represents the drift of the subsidiary's wealth, $\sigma$ represents the uncertainty induced by the hacking on the financial market impacting the portfolio value of the subsidiary with possible accident given by a Poisson random measures $\mu_P$, which intensity $\lambda$ depends on the compromised devices and the direct hacking activity, reflecting the L-hop modeling.
{
 Unlike $\eta$, the Poisson measure $\mu_p$ represents discrete, observable, and sudden catastrophic cyber incidents (e.g., a major ransomware execution or a sudden data breach that triggers an immediate drop in value). While it is mathematically possible to model the arrivals of attacks affecting the SIR dynamics purely through jumps, separating the ambient ambiguity ($\eta$) from the discrete catastrophes ($\mu_p$) allows us to cleanly capture the dual nature of modern cyber threats. 
    }
    Finally, Cyberattackers continuously evolve their tactics, techniques, and procedures. Attackers may exploit vulnerabilities or create novel attack patterns that were not present in the training data, leading to model uncertainty and ambiguity on their actions $\eta$. This issue is usually addressed by adopting a robust approach of the problem; see, for example, \cite{balter2023model,bielecki2014dynamic,vuca,moralhazardambiguity,sung2022optimal}. Let $(\eta, \mathbb{P})$ represent a probability model defined by the cyber attack, leading to the formulation of a Stackelberg bilevel stochastic optimization problem, which can be broadly outlined as follows:
$$\begin{cases}
V_0^P=\sup_{\xi,\hat \alpha}\inf_{(\mathbb P,\eta)} \mathbb E^{\mathbb P}[U_P(\xi,P_T,S_T,I_T,C_T,\hat \alpha,\eta)],\\
\text{ subject to } \\
(IC-\sigma)\quad V_0^A(\xi):=\sup_{ \alpha}\inf_{(\mathbb P,\eta)} \mathbb  E^{\mathbb P}[U_A(\xi,P_T,S_T,I_T,C^A_T,\alpha)]\\
\qquad \qquad \qquad \qquad=\mathbb E^{\mathbb P^{\hat \eta}}[U_A(\xi,P_T,S_T,I_T,C^A_T,\hat \alpha)]\\
(R)\quad V_0^A(\xi)\geq R_0.
\end{cases}$$
We call this problem $\mathbf{(2Mm-\sigma)}$ standing for bilevel Max-min optimization with ambiguity, $(IC-\sigma)$ is the incentive compatibility condition with ambiguity, $(R)$ is the reservation utility constraint, $U_P,U_A$ are the utility functions of the holding company and the subsidiary, respectively, $\xi$ represents the compensation proposed to the subsidiary, and $C_T,C_T^A$ represent the additional discontinuous costs incurred by the holding company and the subsidiary, respectively, as a result of cyber attacks.

\subsection{Comparison with the literature} We now detail the main contributions of this work on three different topics: cyber risk modeling, stochastic optimization and agency problem and cyber risk economics. 
\begin{itemize}
\item \textit{Cyber risk modeling and economics.} While most models studied to date have focused on either discrete-time optimization or deterministic SIR models for cyber risk, our approach addresses cyber risk uncertainty through a fully stochastic framework that includes volatility uncertainty in both the SIR system and the wealth process. This extends, for example, the work of \cite{khouzani2019scalable, hillairet2022cyber}. In addition, we provide a comprehensive model of L-hop risk propagation using a stochastic SIR system with model ambiguity.

Incentive mechanisms for cyber risk management have been previously studied in contexts such as health data protection and optimal cybersecurity investments; see \cite{khouzani2019scalable, zhang2021bayesian, wessels2021understanding, bauer2009cybersecurity, lee2022optimally}. We contribute to this literature by extending the analysis to a continuous-time setting, focusing on the optimal design of incentive schemes using a bilevel max-min optimization approach within a Stackelberg game framework.

\item \textit{Agency problem, stochastic control and optimization.} Stochastic bilevel optimization in continuous time with ambiguity has been previously studied in \cite{sung2022optimal, moralhazardambiguity, vuca}. In this work, we extend this framework to a stochastic bilevel max-min optimization problem in continuous time and volatility uncertainty with jumps. Specifically, we propose a novel connection between second-order backward stochastic differential equations with jumps (2BSDEJs) and principal-agent problems involving both moral hazard and model ambiguity. 2BSDEs have been extensively studied in the literature since the pioneering works \cite{soner2012wellposedness, cheridito2007second, possamai2018stochastic}; see also \cite{popier2019second, possamai2015weak, matoussi2014second}, and more recently, their extensions to include jump processes \cite{Kazi_Tani_2015, denis:hal-04822047,possamai2025mind}. However, the link between 2BSDEs with jumps and principal-agent problems under volatility uncertainty and accident risk has not yet been established. This paper addresses that gap. In particular, we extend the framework of \cite{vuca} to incorporate accidents, and generalize the models in \cite{capponi2015dynamic, bensalem2020prevention} by introducing volatility ambiguity in the context of cyber risk. Finally, we develop a Perron’s method to prove the existence of a viscosity solution to an integro-partial Hamilton–Jacobi–Bellman–Isaacs (HJBI) equation characterized by the principal's value function 
$V_0^P$. This extends the methods in \cite{perronsirbu,bayraktar2012stochastic} and \cite{vuca} to settings with jump-diffusion processes.
\end{itemize}
{
Our technical contributions are twofold:
\begin{enumerate}
    \item Our work directly extend \cite{vuca} to jumps, benefiting from the recent results on 2BSDE with jumps \cite{denis:hal-04822047,possamai2025mind,gennaro2BSDE2025} in non-Markovian framework (driven by the contract $\xi$). This paper is the first connecting a non-Markovian 2BSDE to a stochastic max-min control problem with jumps and ambiguity on the volatility, see Theorem \ref{thm:agent-2BSDEJ};
\item as far as we know this paper, and its first version is the first extending \cite{cvitanic2018dynamic,vuca} to a robust study of principal-agent problem with jumps and volatility ambiguity and making a connection between the agents problem and the solution to a 2BSDE with jumps (Theorem 3.7), reducing the principal problem to the solution of an Isaacs-HJB equation (fully non linear integro-PDE and using Perron's method (see Appendix \ref{app:viscosity_proof}) with jumps in the context of max-min to characterize the viscosity property of the solution to the problem, extending \cite{bayraktar2012stochastic,BayraktarSirbu2013}, see Theorem \ref{thm:viscosity}.
\end{enumerate}
}
The structure of this work is as follows. Section \ref{section:model} presents the modeling framework, including the canonical process and weak formulation of the problem, the controlled equation, admissible controls and contracts, and finally the bilevel max-min stochastic optimization. Section \ref{section:subsidiary} focuses on the incentive compatibility (IC) condition, also known as the agent's problem and its connection to a 2BSDE with jumps. Section \ref{section:holding} investigates the optimal compensation schemes by reducing the problem to an integro-Isaacs PDE, applying a verification theorem and Perron’s method in the context of discontinuous stochastic processes. Finally Section \ref{sec:cyber} applies the results to cyber risk management illustrated with numerical experiments exploring the benefit of a contracting mechanism to monitor both the cyber threat and its uncertainty. 
\newpage
\section{The model and bilevel max-min problem}
\label{section:model}
\subsection{Canonical process and weak formulation}\label{section:modgen}

Fix a horizon $T>0$ and integers
$$
n\ \text{(state dimension)},\qquad \ell\ \text{(Brownian dimension)},\qquad m\ \text{(mark dimension)}.
$$
Let $E\subset\mathbb{R}^m\setminus\{0\}$ be a Borel mark space with Borel $\sigma$-algebra $\mathcal B(E)$.
We fix a predictable \emph{base compensator} $\nu^0_t(de)dt$ on $[0,T]\times E$, which is $\sigma$-finite and has full support on $E$. Define
$$
\Omega^c := \big\{\omega\in C([0,T];\mathbb{R}^{\ell}) : \omega_0=0\big\},\quad
\Omega^d := \mathsf M_p\big((0,T]\times E\big),\quad
\Omega^x := D([0,T];\mathbb{R}^n),
$$
where $\mathsf M_p((0,T]\times E)$ is the space of \emph{integer–valued measures} on $(0,T]\times E)$.
We equip $\Omega^c$ with the Wiener topology, $\Omega^d$ with the vague topology, and $\Omega^x$ with the Skorokhod topology; in particular, each factor is Polish.
Set the canonical product space and its $\sigma$-field
$$
\Omega := \Omega^c \times \Omega^d \times \Omega^x,\qquad
\mathcal G := \mathcal B(\Omega^c)\otimes\mathcal B(\Omega^d)\otimes\mathcal B(\Omega^x).
$$
On $(\Omega,\mathcal G)$, define the coordinate processes
$$
W^0_t(\omega^c,\omega^d,\omega^x) := \omega^c_t,\qquad
\mu^0(B)(\omega^c,\omega^d,\omega^x) := \omega^d(B),\quad B\in\mathcal B\big((0,T]\times E\big),
$$
$$
X_t(\omega^c,\omega^d,\omega^x) := \omega^x_t,\qquad t\in[0,T].
$$
Thus $W^0$ is the Brownian coordinate, $\mu^0$ the jump (integer–valued) coordinate with base compensator $\nu^0_t(de)dt$, and $X$ the state coordinate. Let $\mathbb G=(\mathcal G_t)_{t\in[0,T]}$ be the raw filtration generated by $(W^0,\mu^0,X)$, i.e.
$$
\mathcal G_t := \sigma\big(W^0_s,\ \mu^0\big((0,s]\times A\big),\ X_s:\ 0\le s\le t,\ A\in\mathcal B(E)\big).
$$
Let $\mathbb G^+=(\mathcal G_t^+)_{t\in[0,T]}$ denote its right–continuous modification, $\mathcal G_t^+ := \bigcap_{u>t}\mathcal G_u$. Let $\mathcal M(\Omega)$ denote the set of all probability measures on $(\Omega,\mathcal G)$.
Define the universal filtration
$$
\mathcal G_t^* := \bigcap_{\mathbb P\in\mathcal M(\Omega)} \mathcal G_t^{\mathbb P},\qquad
\mathbb G^* := (\mathcal G_t^*)_{t\in[0,T]},
$$
where $\mathcal G_t^{\mathbb P}$ is the usual augmentation of $\mathcal G_t$ under $\mathbb P$. For $\mathbb P\in\mathcal M(\Omega):\quad
\mathbb F^{\mathbb P} := (\mathcal F_t^{\mathbb P})_{t\in[0,T]}$ is the right–continuous, $\mathbb P$–complete augmentation of $\mathbb G$. For a nonempty $\mathcal P\subset\mathcal M(\Omega)$, a set $N\in\mathcal G$ is \emph{$\mathcal P$-polar} if $\mathbb P(N)=0$ for all $\mathbb P\in\mathcal P$.
Let $\mathcal T^{\mathcal P}$ be the $\sigma$-algebra of $\mathcal P$-polar sets and define the $\mathcal P$–universal filtration
$$
\mathcal F_t^{\mathcal P} := \mathcal G_t^* \vee \mathcal T^{\mathcal P},\qquad
\mathbb F^{\mathcal P} := (\mathcal F_t^{\mathcal P})_{t\in[0,T]},
$$
with right–continuous modification $\mathbb F^{\mathcal P,+}$.
When harmless, we omit the superscript $\mathcal P$. For $\mathcal P\subset\mathcal M(\Omega)$, $t\in[0,T]$, and $\mathbb P\in\mathcal P$, set
$$
\mathcal P[\mathbb P,\mathbb F^+,t] := \big\{\mathbb P'\in\mathcal P:\ \mathbb P'=\mathbb P\ \text{on}\ \mathcal F_t^+\big\}.
$$
For any $\mathbb P\in\mathcal M(\Omega)$ and any $\mathbb F$–stopping time $\tau$, there exists a family of regular conditional probabilities $(\mathbb P^\tau_\omega)_{\omega\in\Omega}$ (standard).

\medskip
It is well known (see, e.g., \cite{stroock1997multidimensional}) that for every $\mathbb{P} \in \mathcal{M}(\Omega)$ and every $\mathbb{F}$-stopping time $\tau$ with values in $[0,T]$, there exists a family of \emph{regular conditional probability distributions (r.c.p.d.)} $(\mathbb{P}_\omega^\tau)_{\omega \in \Omega}$; we refer to \cite[Section~1.1.3]{possamai2018stochastic} for details.

\begin{definition}[Admissible laws with fixed jump law]\label{def:PS-fixed-jump}
Fix a predictable base compensator $\nu^0_t(de)dt$ on $[0,T]\times E$. Let
$$
\lambda^0_t(d\chi):=\int_E \mathbf 1_{\{\beta(t,X_{t-},e)\in d\chi\}}\ \nu^0_t(de), \qquad \Lambda^0(ds,d\chi):=\lambda^0_t(d\chi)ds.
$$
be the state–dependent compensator on $\mathbb R^n\setminus\{0\}$ induced by $\beta$. For $t\in[0,T]$ and $x\in\mathbb R^n$, define $\mathcal P(t,x)$ as the set of $\mathbb P\in\mathcal M(\Omega)$ such that:

\begin{enumerate}
\item[(i)] Under $\mathbb P$, $W^0$ is an $\ell$-dimensional $\mathbb F^\mathbb P$–Brownian motion, $\mu^0$ is integer–valued with predictable compensator $\nu^0_t(de)dt$, and $W^0$ is independent of $\mu^0$.

\item[(ii)]
Under $\mathbb P$, $X$ is an $\mathbb F^{\mathbb P}$–semimartingale with canonical decomposition
$$
X_t
= X_0 + X^{c,\mathbb P}_t
+ \int_{(0,t]\times(\mathbb R^n\setminus\{0\})} \chi\big(\mu_X-\Lambda^{0}\big)(ds,d\chi),
\qquad t\in[0,T],
$$
where $\mu_X$ is the jump measure of $X$ and $\lambda^{0}$ is its $\mathbb F^{\mathbb P}$–predictable compensator defined above, satisfying
$$
\int_0^T\int_{\mathbb R^n\setminus\{0\}} (1\wedge |\chi|^2)\Lambda^{0}(ds,d\chi)<\infty,
\qquad \mathbb P\text{-a.s.}
$$

\item[(iii)] 
$\langle X^{c,\mathbb P}\rangle_t = \int_0^t \widehat\sigma^{\mathbb P}_sds$ for a predictable
$\widehat\sigma^{\mathbb P}\in\mathbb S_+^n$.
\end{enumerate}

Specifically, we denote $\mathcal P:=\mathcal P(0,x_0)$.
\end{definition}

We then have the following lemma, whose proof follows the same line as in the proof of \cite[Proposition 5.3]{cvitanic2018dynamic}.

\begin{lemma}
    By construction, $\mathcal P(t,x)$ is \emph{saturated}: if $\mathbb P\in\mathcal P(t,x)$ and $\mathbb Q\sim\mathbb P$ under which $X$ is a local martingale, then $\mathbb Q\in\mathcal P(t,x)$.
\end{lemma}

It is well known (see, e.g., \cite{KARANDIKAR199511}) there exists an $\mathbb F$–progressively measurable aggregator $\langle X\rangle$ whose continuous density
$$
\widehat\sigma_t:=\limsup_{\varepsilon\downarrow 0}\frac{\langle X\rangle^c_t-\langle X\rangle^c_{t-\varepsilon}}{\varepsilon}\in\mathbb S_+^n
$$
satisfies $\widehat\sigma_t=\widehat\sigma^{\mathbb P}_t$ for $dt\otimes d\mathbb P$–a.e.\ $(t,\omega)$ and all $\mathbb P\in\mathcal P$.

\subsection{Admissible controls and Girsanov via Doléans–Dade exponentials}

\begin{assumption}[Regularity on model's data]\label{ass:standing-data}
Fix compact metric spaces $A$ and $H$. Let
$$
b:[0,T]\times\Omega\times A\times H\to\mathbb R^n,\quad
\sigma:[0,T]\times\Omega\times H\to \mathcal M_{n,\ell}(\mathbb R),\quad
\beta:[0,T]\times\Omega\times E\to\mathbb R^n,
$$
be $\mathbb F$–predictable in $(t,\omega)$ and continuous in the control arguments. Set
$
\Sigma(t,\omega,h):=\sigma\sigma^\top(t,\omega,h)\in\mathbb S_+^n.
$ We assume
\begin{enumerate}
\item (Growth/Lipschitz) $b,\sigma$ are locally bounded and Lipschitz in the state, uniformly on compact control sets. There exists $0<\bar\kappa$ such that $$\bigl\| {b}(t, x, a, h) \bigr\|\le
\bar\kappa\Bigl(1 + \|x\|_{t,\infty} + |a| \Bigr),
\quad
\bigl\|\partial_{ a}b(t, x, a, h)\bigr\|
\le
\bar\kappa\Bigr).$$
\item (Jump integrability) $\displaystyle \int_E(1\wedge|\beta(t,\omega,e)|^2)\nu^0_t(de)<\infty$ for all $(t,\omega)$.
\item (Base compensator) $\nu^0_t(de)dt$ is a fixed predictable compensator on $[0,T]\times E$ with full support on $E$.
\item (Covariance realization) There exists an $\mathbb F$–predictable process $\eta$ with values in $H$ such that
$$
\Sigma\big(t,X_{t-},\eta_t\big)=\widehat\sigma_t
\quad\text{for }dt\otimes d\mathcal P\text{-q.s.}
$$
\end{enumerate}
\end{assumption}

{
\begin{remark}
The requirement 4. that the quadratic variation density $\hat{\sigma}_t$ is realized by a process $\eta_t$ carries a fundamental economic meaning in our framework. Mathematically,  $\eta$ acts as a parametrization of the volatility uncertainty. By imposing this structure, we are effectively restricting the set of admissible priors (the family of probability measures considered for the weak formulation) to those whose volatility matrices take the form $\Sigma(t, x, h)$ for some $h \in H$. Economically, in the context of our max-min bilevel problem, this mathematical parametrization coincides with the control of an adversary. Rather than just being an abstract label for the ``worst-case volatility'' chosen by a fictitious Nature, $\eta$ models a concrete adversarial action. For instance, in the cyber risk application developed in Section \ref{sec:cyber}, $\eta$ naturally represents the attack effort of a hacker aiming to maximize the system's vulnerability.

\end{remark}
}

As usual in moral hazard contract theory, see \cite{cvitanic2017moral,mastrolia2025agency} the agent modifies the distribution of the canonical process by changing the reference probability measure $\mathbb P^{0}\in \mathcal P(0,x_0)$ to a new probability measure $\mathbb P^{\alpha,\eta}$. We then define the set of admissible controls and feasible priors through the Girsanov Theorem.

\begin{definition}[Admissible controls]\label{def:admissible-controls}
A pair $(\alpha,\eta)$ of $\mathbb F$–predictable processes with values in $A\times H$ is \emph{admissible} if there exist predictable processes $\kappa(e;\alpha,\eta)$ and $\zeta(\alpha,\eta)$ such that
$\kappa_t(e;\alpha,\eta) > 0$ with\footnote{${\Sigma}^\dagger$ denote the Moore–Penrose pseudoinverse. Specificaly, if $\sigma$ has full row rank, then it is $\Sigma^{-1}$; if $\sigma$ has full column rank, then ${\Sigma}^\dagger = (\sigma^\dagger)^\top\sigma^\dagger = \sigma(\sigma^\top\sigma)^{-2}\sigma^\top$}
   $$
    \int_0^T \int_E \big(\sqrt{\kappa_t(e;\alpha,\eta)} - 1\big)^2 \nu^0_t(de)dt < \infty,
    $$
\begin{equation}
    \label{eq:zeta-solution}
    \zeta_t(\alpha,\eta)
    = \Sigma^\dagger(t,X_{t-},\eta_t)\left(b(t,X_{t-},\alpha_t,\eta_t)
    - \int_E \beta(t,X_{t-},e)\big(\kappa_t(e;\alpha,\eta) - 1\big)\nu^0_t(de)\right),
    \end{equation}
and
    $
b(t,X_{t-},\alpha_t,\eta_t)
    - \int_E \beta(t,X_{t-},e)\big(\kappa_t(e;\alpha,\eta) - 1\big)\nu^0_t(de)\in\mathrm{Ran}\Sigma
$
where 
$$
    b(t,X_{t-},\alpha_t,\eta_t)
    = \Sigma(t,X_{t-},\eta_t)\zeta_t(\alpha,\eta)
    + \int_E \beta(t,X_{t-},e)\big(\kappa_t(e;\alpha,\eta) - 1\big)\nu^0_t(de),
$$
and such that 
\begin{equation}\label{eq:expo}
\mathbb E^{\mathbb P^{0}}\Big[\mathcal E\left(\int_0^\cdot\int_E(\kappa-1)\tilde\mu^0(ds,de)\right)_T
 \mathcal E\left(\int_0^\cdot \zeta_s(\alpha,\eta)^\top dX^{c}_s\right)_T\Big]=1,
 \end{equation}
where $\mathcal E()$ denotes the Doleans-Dade exponential process: 
\begin{align*}
\mathcal E&\Big(\int_0^\cdot\int_E(\kappa-1)\tilde\mu^0(ds,de)\Big)_t\\
:= \exp&\left(
\int_0^t\int_E \log \kappa_s(e)\mu^0(ds,de)
- \int_0^t\int_E \big(\kappa_s(e)-1\big)\nu^0_s(de)ds
\right),
\end{align*}
and
$$
\mathcal E\Big(\int_0^\cdot \zeta_s(\alpha,\eta)^\top dX^{c}_s\Big)_t
= \exp\left(
\int_0^t \zeta_s(\alpha,\eta)^\top dX^{c}_s
- \tfrac12 \int_0^t \zeta_s(\alpha,\eta)^\top  d[X^c]_s  \zeta_s(\alpha,\eta)
\right),
$$ and we denote the compensated measure as $\tilde{\mu}^0(dt,de):=\mu^0(dt,de)-\nu^0_t(de)dt$.
\end{definition}
As a consequence of the admissibility of $\alpha,\eta$ we can define 
\begin{equation}
    \label{eq:nu-solution}
    \nu_t^{(\alpha,\eta)}(de) := \kappa_t(e;\alpha,\eta)\nu^0_t(de),
    \end{equation}
    and a probability $\mathbb P^{(\alpha,\eta)}$ by
$$
\frac{d\mathbb P^{(\alpha,\eta)}}{d\mathbb P^{0}}
=\mathcal E\left(\int_0^\cdot\int_E(\kappa-1)\tilde\mu^0(ds,de)\right)_T
 \mathcal E\left(\int_0^\cdot \zeta_s(\alpha,\eta)^\top dX^{c}_s\right)_T,
$$
under which
$$
dX_t=b(t,X_{t-},\alpha_t,\eta_t)dt
   +\sigma(t,X_{t-},\eta_t)dW^{(\alpha,\eta)}_t
   +\int_E \beta(t,X_{t-},e)\tilde\mu^{(\alpha,\eta)}(dt,de),
$$
with jump compensator $\nu^{(\alpha,\eta)}_t=\kappa_t(\cdot;\alpha,\eta)\nu^0_t$ and corresponding compensated measure $\tilde\mu^{(\alpha,\eta)}(dt,de):=\mu^{(\alpha,\eta)}(dt,de)-\nu_t^{(\alpha,\eta)}(de)dt$. { For each $\alpha$, we define $$\mathcal H^{(\alpha)}:=\{(\mathbb P^{(\alpha,\eta)},\eta):\eta\in\mathfrak H\}, \quad\mathcal P^{(\alpha)} := \{\mathbb P^{(\alpha,\eta)}:\eta\in\mathfrak H\}.$$}{\begin{remark}[Conditions ensuring \eqref{eq:expo}]Recall that \eqref{eq:expo} requires the Dol\'eans--Dade exponential defining the change of
measure to be a true $(\mathbb P_0,\mathbb F)$--martingale.
For example, we can assume that there exist constants $\bar\zeta>0$,
$0<\underline\kappa\le \overline\kappa<\infty$, and $C<\infty$ such that, for every admissible
pair $(\alpha,\eta)$,
\begin{align*}
&|\zeta_t(\alpha,\eta)| \le \bar\zeta,
\qquad
\underline\kappa \le \kappa_t(e;\alpha,\eta)\le \overline\kappa,
\quad dt\otimes d\mathbb P_0\otimes \nu_t^0(de)\text{-a.e.},\\
&\int_0^T \zeta_s(\alpha,\eta)^\top \Sigma(s,X_{s-},\eta_s)\,\zeta_s(\alpha,\eta)\,ds \le C,
\qquad \mathbb P_0\text{-a.s.},\\
&\int_0^T\int_E \Big(\kappa_s(e;\alpha,\eta)\log \kappa_s(e;\alpha,\eta)
-\kappa_s(e;\alpha,\eta)+1\Big)\,\nu_s^0(de)\,ds \le C,
\qquad \mathbb P_0\text{-a.s.}
\end{align*}

Under these conditions, the continuous exponential is a true martingale by a Novikov/Kazamaki
criterion, and the jump exponential is a true martingale by a L\'epingle--M\'emin type criterion, see \cite{lepingle1978integrability,okada1982criterion}, \cite[Theorem 1.31]{oksendal2019applied},
for purely discontinuous local martingales. Consequently, (2.2)
is satisfied.

\end{remark}}

\begin{remark}
If ${\Sigma}(t,x,\eta)$ is uniformly elliptic, then $\zeta(\alpha,\eta)$ in \eqref{eq:zeta-solution} is unique and given by the usual inverse; otherwise the range condition above is the natural compatibility restriction for attainable drifts. Uniform bounds and compactness of $A,H$, together with continuity of coefficients, imply compactness of the attainable covariance set and ensure the Novikov/Lépingle–Mémin criteria can be enforced uniformly.
\end{remark}

\begin{remark}
    In the classical framework, as in \cite{moralhazardambiguity,vuca}, the Principal and Agent may hold different beliefs about the volatility, leading to distinct sets of admissible laws. However, in our problem setup, particularly in the context of a holding company and its subsidiary, it is customary to assume that they share the same belief.
\end{remark}

\subsection{Bi-level optimization: agent best response and principal’s problem with volatility \& jump control}

The principal offers an $\mathcal F_T$-measurable compensation $\xi$.
Let the discount factor be
$$
\mathcal K_{t,s}:=\exp\left(-\int_t^s k(r,X_r)dr\right),\qquad 0\le t\le s\le T,
$$
for a given predictable rate $k$.
We assume $\xi$ belongs to
$$
\Xi:=\Big\{\xi\in L^0(\mathcal F_T):\ 
\sup_{\mathbb P\in\mathcal P}\ 
\mathbb E^{\mathbb P}\big[\mathcal K_{0,T}\big(|U^A(\xi)|+|F^A(X_T)|\big)\big]<\infty\Big\},
$$
{where $U^A$ is a nondecreasing and concave utility function and $F^A$ has polynomial growth in $X$.} Given a contract $\xi\in\Xi$, the agent’s worst-case value is
\begin{equation}\label{eq:agent-value}
{\
\begin{aligned}
V_0^A(\xi)
:=\ \sup_{\alpha\in\mathfrak A}\ \inf_{(\mathbb P,\eta)\in \mathcal H^{(\alpha)}}
\ \mathbb E^{\mathbb P}\Big[
&\ \mathcal K_{0,T}\big(U^A(\xi)+F^A(X_T)\big)\\
&\ -\ \int_0^T \mathcal K_{0,s} C^A\big(s,X_s,\alpha_s,\eta_s\big)ds
\Big].
\end{aligned}}
\end{equation}
A measurable selection $\widehat\alpha(\xi)\in\mathfrak A$ with\footnote{{
    Since we assume $A$ is compact (Assumption~2.3) and the agent Hamiltonian and the objective is $\mathcal{P}\otimes\mathcal{B}(A)$--measurable in $(t,\omega,a)$ and continuous in $a$, the argmax correspondence is nonempty and compact-valued and admits a predictable measurable maximizer by a measurable maximum theorem, see Theorem 18.19~\cite{alma991010524309706532} and Theorem 2~\cite{Brown1973Measurable} and more generally \cite{karoui2013capacitiesI,karoui2013capacities}.

}}
$$
\widehat\alpha(\xi)\in\arg\max_{\alpha\in\mathfrak A}\ \inf_{(\mathbb P,\eta)\in \mathcal H^{(\alpha)}}
\ \mathbb E^{\mathbb P}\Big[
\mathcal K_{0,T}\big(U^A(\xi)+F^A(X_T)\big)
-\int_0^T \mathcal K_{0,s} C^A\big(s,X_s,\alpha_s,\eta_s\big)ds\Big]
$$
is called an agent best response. The individual rationality constraint is
\begin{equation}\label{eq:IR}
{\qquad V_0^A(\xi)\ \ge\ R_0,\qquad}
\end{equation}
for a given reservation level $R_0$.

\begin{remark}
If one models additional running costs via marked Poisson processes (e.g., $N^A,N^P$),
then under linear expectation and dominated jumps those costs can indeed be absorbed into
$C^A,C^P$ by taking expectations:
$$
\mathbb E\left[\int_0^T\int L(t,X_t,\cdot)N(dt,de)\right]
=\mathbb E\left[\int_0^T\int L(t,X_t,\cdot)\lambda(t,X_t,\cdot)\nu(de)dt\right].
$$
Thus writing the aggregated forms $C^A,C^P$ is without loss for the problems
\eqref{eq:agent-value}.
\end{remark}

{Let $F^P:\mathbb R^n\to\mathbb R$ be the principal's terminal payoff. 
We assume that $F^P$ is continuous and has polynomial growth, i.e., there exist
constants $C>0$ and $m\ge 1$ such that
$$
    |F^P(x)|\le C(1+|x|^m),\qquad x\in\mathbb R^n .
$$}
The principal chooses $\xi\in\Xi$ to maximize her worst-case expected utility given the agent’s best response:
\begin{equation}\label{eq:principal-problem}
\begin{aligned}
V_0^P
:=&\sup_{\xi\in\Xi}\sup_{\hat\alpha(\xi)} \inf_{(\mathbb P,\eta)\in \mathcal H^{(\widehat\alpha(\xi))}}
\mathbb E^{\mathbb P}\Big[
F^P(X_T)-\xi-\int_0^T C^P\big(s,X_s,\widehat\alpha_s(\xi),\eta_s\big)ds
\Big]
\\
\text{subject to }\quad
&\text{(IC)}\ \ \widehat\alpha(\xi)\ \text{is an agent best response in \eqref{eq:agent-value}},\\
&\text{(IR)}\ \ V_0^A(\xi)\ \ge R_0\ \text{ as in \eqref{eq:IR}}.
\end{aligned}
\end{equation}

\section{Solving the agent’s problem via 2BSDE with jumps}
\label{section:subsidiary}

\subsection{Agent driver and covariance-constrained Hamiltonians}
\label{subsec:agent-driver}
For $$(t,x,y,z,u,a,h) \in [0,T]\times\mathbb R^n\times\mathbb R\times\mathbb R^n
\times \mathcal L_\nu^p\times A\times H,$$ we set
\begin{equation}\label{eq:G}
\begin{split}
G(t,x,y,z,u;a,h)
:= &-k(t,x)y - C^A(t,x,a,h)+ b(t,x,a,h)\cdot z\\
&{+
\int_E
\Bigl[
u(\beta(t,x,e))
-
z\cdot\beta(t,x,e)
\Bigr]
\bigl(\kappa_t(e;a,h)-1\bigr)\nu_t^0(de)}.
\end{split}
\end{equation}

For $(t,x)\in[0,T]\times\mathbb R^n$ and $\Sigma\in\mathcal S_n^+$, define
$$
\mathcal H(t,x,\Sigma):=\{h\in H:\ \sigma{\sigma}^\top(t,x,h)=\Sigma\},
$$ and denote by $\mathcal H(\hat\sigma)$ the set of control $\eta\in \mathfrak H$ with values in $\mathcal H(t,x,\hat\sigma)$, $dt\otimes \mathbb P$ a.e., for every $\mathbb P\in \mathcal P$. We also define the optimized driver at fixed covariance as
\begin{equation}\label{eq:Gstar}
{
G^*(t,x,y,z,u;\Sigma)
:=\ \sup_{a\in A}\ \inf_{h\in \mathcal H(t,x,\Sigma)}\ G(t,x,y,z,u;a,h).
}
\end{equation}

\begin{assumption}[Isaacs at fixed covariance]\label{ass:Isaacs-fixed-Sigma}
For all $(t,x,y,z,u,\Sigma)$,
$$
\inf_{h\in \mathcal H(t,x,\Sigma)}\ \sup_{a\in A} G(t,x,y,z,u;\,a,h)
=\sup_{a\in A}\ \inf_{h\in \mathcal H(t,x,\Sigma)} G(t,x,y,z,u;\,a,h).
$$
\end{assumption}

{ 
To assure the Issacs condition \ref{ass:Isaacs-fixed-Sigma} holds, we provide sufficient conditions and examples below.
\begin{example}
Fix $(t,x,y,z,u,\Sigma)$. Assume $A\subset\mathbb{R}^{d_a}$ and $\mathcal H(t,x,\Sigma)\subset\mathbb{R}^{d_h}$
are nonempty, compact, and convex. Suppose $(a,h)\mapsto G(t,x,y,z,u;a,h)$ defined in (3.1)
is jointly continuous on $A\times \mathcal H(t,x,\Sigma)$ and satisfies:
\begin{enumerate}
\item for every $h\in \mathcal H(t,x,\Sigma)$, the map
$a\mapsto G(t,x,y,z,u;a,h)$ is concave on $A$.
\item for every $a\in A$, the map
$h\mapsto G(t,x,y,z,u;a,h)$ is convex on $\mathcal H(t,x,\Sigma)$.
\end{enumerate}
Then the Isaacs equality holds:
$$
\inf_{h\in \mathcal H(t,x,\Sigma)}\sup_{a\in A}G(t,x,y,z,u;a,h)
=
\sup_{a\in A}\inf_{h\in \mathcal H(t,x,\Sigma)}G(t,x,y,z,u;a,h).
$$
\end{example}

\begin{example}
\label{ex:lq_jumps_isaacs}
Fix $(t,x,y,z,u,\Sigma)$ and suppose $A\subset\mathbb{R}^{d_a}$ and $H(t,x,\Sigma)\subset\mathbb{R}^{d_h}$
are compact convex.
Assume the local part of the driver takes the linear--quadratic form
$$
G_{\mathrm{loc}}(a,h)
=
\langle \theta_a, a\rangle - \tfrac12 a^\top R a
\;+\;
\langle \theta_h, h\rangle + \tfrac12 h^\top Q h
\;+\; \text{(terms independent of $(a,h)$)},
$$
where $R\in\mathbb{S}^{d_a}_{++}$ and $Q\in\mathbb{S}^{d_h}_{++}$, and where
$\theta_a,\theta_h$ may depend on $(t,x,y,z,\Sigma)$ (e.g.\ through the drift term $b\cdot z$). This is for example satisfied if $C^A$ is quadratically separable and $b$ is linearly separable in $a,h$. For the jump contribution, assume the compensator density is affine in $(a,h)$, so the jump term
$$
\int_E
\Bigl[
u(\beta(t,x,e))
-
z\cdot\beta(t,x,e)
\Bigr]
\bigl(\kappa_t(e;a,h)-1\bigr)\nu_t^0(de)
$$
is affine in $(a,h)$ and contributes no curvature. Consequently, $a\mapsto G(a,h)=G_{\mathrm{loc}}(a,h)+G_{\mathrm{jump}}(a,h)$ is concave
and $h\mapsto G(a,h)$ is convex, and (Sufficient condition) yields Isaacs equality,
so Assumption~3.1 holds. Moreover, in the unconstrained case the saddle point is explicit:
$a^\ast = R^{-1}\theta_a,\, h^\ast = -Q^{-1}\theta_h,$
and under compact constraints it is given by projection onto $A$ and $H(t,x,\Sigma)$, respectively.
\end{example}}

We then define the Hamiltonian $H: [0, T] \times \mathbb R^n \times \mathbb{R} \times \mathbb{R}^n\times  \mathcal L^{p,m}_\nu \times  \mathcal S_n^+ \to \mathbb{R}$
\begin{equation}\label{eq:H}
H(t,x,y,z,u;\Gamma)\ :=\ \inf_{\Sigma\in\mathbb S_+^n}\Big\{\tfrac12\mathrm{Tr}(\Sigma\Gamma)\ +\ G^*(t,x,y,z,u;\Sigma)\Big\}.
\end{equation}

\subsection{2BSDE with jumps for the agent and verification}
\label{subsec:2bsdej-agent}

Given $\xi\in\Xi$, the 2BSDEJ reads
\begin{equation}\label{eq:2BSDEJ-agent}
\begin{aligned}
Y_t &= U^A(\xi)+F^A(X_T)
+\int_t^T G^{*}\big(s,X_s,Y_s,Z_s,U_s;\widehat\sigma_s\big)ds \\
&\quad -\int_t^T Z_s\cdot dX^{c,\mathbb P}_s
      -\int_t^T\int_{\mathbb R^n\setminus\{0\}} U_s(\chi)(\mu_X-\Lambda^0)(ds,d\chi)
      -\int_t^T dK_s^{\mathbb P},
\quad \mathcal P - q.s..
\end{aligned}
\end{equation}

\begin{definition}\label{2bsde:def}
    We say that a quadruplet $(Y, Z, U, K)$ is a solution to the 2BSDEJ \eqref{eq:2BSDEJ-agent} if there exists $p > 1$ such that 
$$
(Y, Z, U, K) \in \mathbb S^p_0(\mathbb{F}^{\mathcal{P}}_+ , \mathcal{P}) \times  \mathbb H^p_0(\mathbb{F}^{\mathcal{P}} , \mathcal{P}) \times \mathbb J^p_0(\mathbb{F}^{\mathcal{P}} , \mathcal{P}) \times \mathbb K^p_0(\mathbb{F}^{\mathcal{P}} , \mathcal{P})
$$
satisfies \eqref{eq:2BSDEJ-agent} and $K$ satisfies the minimality condition
\begin{equation}
\label{eq:min.con}
0= \operatorname*{essinf}_{\mathbb P' \in \mathcal{P}[\mathbb P,\mathbb{F}+,s]} \mathbb{E}^{\mathbb P'}\left[ K_T-K_s \mid \mathcal{F}_{s}^{\mathbb P,+} \right], \quad s \in [t, T], \ \mathbb P \text{ - a.s.}, \ \forall \mathbb P \in \mathcal{P}.
\end{equation}
\end{definition} 

\begin{assumption}[Regularity for well–posedness]\label{ass:agent-regularity}
$A,H$ are compact; $b,k,C^A$ are bounded on compacts and continuous in $(a,h)$; $\sigma$ is bounded and continuous in $(a,h)$ so that the attainable covariance correspondence has compact values; moreover, $G^*$ is locally Lipschitz in $(y,z)$, uniformly on compacts.
\end{assumption}

\begin{lemma}\label{lem:2bsdej-wp}
Under Assumptions~\ref{ass:Isaacs-fixed-Sigma} and \ref{ass:agent-regularity}, for any $\xi\in\Xi$ with $U^A(\xi)\in\mathbb L^{p,\kappa}_0$, the 2BSDEJ \eqref{eq:2BSDEJ-agent} admits a unique solution
$$
(Y,Z,U,K)\in
\mathbb S^p_0(\mathbb F^{\mathcal P,+},\mathcal P)\times
\mathbb H^p_0(\mathbb F^{\mathcal P},\mathcal P)\times
\mathbb J^p_0(\mathbb F^{\mathcal P},\mathcal P)\times
\mathbb K^p_0(\mathbb F^{\mathcal P},\mathcal P).
$$
\end{lemma}

\begin{theorem}\label{thm:agent-2BSDEJ}
Let $(Y,Z,U,K)$ solve \eqref{eq:2BSDEJ-agent}.
Then $Y_0$ is $\mathcal F_0$–measurable and constant under every $\mathbb P\in\mathcal P(0,x_0)$, and
\begin{equation}
\label{eq:ver}
V_0^A(\xi)=\sup_{a\in A}\inf_{(\mathbb P,\eta)\in\mathcal H^{(a)}}\mathbb E[Y_0].
\end{equation}
Moreover, a triplet $(\widehat\alpha,\widehat\eta,\widehat{\mathbb P})$ is optimal if and only if
$$
G^*\big(t,X_t,Y_t,Z_t,U_t;\widehat{\sigma}_t\big)
=G\big(t,X_t,Y_t,Z_t,U_t;\widehat\alpha_t,\widehat\eta_t\big)
\quad\text{for }dt\otimes d\widehat{\mathbb P}\text{-a.e.,}
\quad
K_T^{\widehat{\mathbb P}}=0\ \ \widehat{\mathbb P}\text{-a.s.}
$$
\end{theorem}

{
\begin{remark}[Interpretation of the process $K$]
Note that the solution of the 2BSDEJ \eqref{eq:2BSDEJ-agent} contains a family of nondecreasing processes $(K^{\mathbb P})_{\mathbb P \in \mathcal{P}}$ indexed by the set of admissible priors. Because the 2BSDE solution $Y$ must satisfy the backward equation $\mathcal{P}$-almost-surely under all measures $\mathbb P \in \mathcal{P}$ simultaneously, $Y$ effectively acts as an upper envelope over all possible scenarios, since $Y$ is a supremum over all priors. When we evaluate the system under the true worst-case measure $\hat{\mathbb P}$ (the optimal control of the adversary), the robust value of the problem matches the exact value, meaning no ``slack'' is needed to satisfy the backward equation, leading to $K^{\hat {\mathbb P}}_T = 0$. In the context of our agency problem, this implies that when the cyber threat inflicts maximum damage, the principal's robust valuation requires no upward correction. Mathematically, 
$$Y_t=\esssup_\alpha\essinf_{\mathbb P,\eta}\;  \mathcal Y_t^{\alpha,\mathbb P,\eta},$$
where $\mathcal Y^{\alpha,\mathbb P,\eta}$ is the first component of the solution to a BSDE with generator $G(\dots,\alpha,\eta)$ and terminal condition $U^A(\xi)$. From Isaac-conditions and classical comparison results for BSDE with jumps, we get 
$$Y_t=  \mathcal Y_t^{\hat\alpha,\hat{\mathbb P},\hat\eta},\; \hat{\mathbb P}-a.s.$$ Consequently, we automatically have $K^{\hat{\mathbb P}}=0,\; \hat{\mathbb P}-a.s.$ by identification and uniqueness of the solution to the 2BSDEJ. 
\end{remark}

Finally, note that the maximizer $\hat\alpha$ may not be unique. We denote by $\hat{\mathcal A}$ the set of optimal $\hat\alpha$ in the sense of the previous theorem. Going back to the initial Principal-Agent problem with moral hazard as formalized by Holmstrom and Milgrom, the Principal controls both the contract and a recommended effort to the agent, choosing the best incentive compatible efforts for her own payoff, see \cite[Equation (7)]{holmstrom1987aggregation}.

}

\section{Optimal contract, Perron's method, and viscosity characterization}
\label{section:holding}

{  Let $(Y,Z,U,K)$ be the solution of 2BSDEJ~\eqref{eq:2BSDEJ-agent}, then \begin{equation}
\begin{split}
    U^A(\xi) + F^A(X_T) = Y_0-\int_t^T G^{*}\big(s,X_s,Y_s,Z_s,U_s;\widehat\sigma_s\big)ds  +\int_t^T Z_s\cdot dX^{c,\mathbb P}_s\\
    +\int_t^T\int_{\mathbb R^n\setminus\{0\}} U_s(\chi)(\mu_X-\Lambda^0)(ds,d\chi)
      +\int_t^T dK_s^{\mathbb P},
\quad \mathcal P - q.s.,
\end{split}
\end{equation} where the process $K$ satisfies the minimality condition \eqref{eq:min.con}. In addition, the reserved utility requirement \eqref{eq:IR}, \textit{i.e.} \begin{equation*}
    \sup_{a\in A}\inf_{(\mathbb P,\eta)\in\mathcal H^{(a)}}\mathbb E[Y_0] \geq R_0,
\end{equation*} is satisfied. We define the set of those $\mathcal F_0-$measurable random variable as \begin{equation*}
    \mathbb Y_0:=\left\{Y_0,\quad \sup_{a\in A}\inf_{(\mathbb P,\eta)\in\mathcal H^{(a)}}\mathbb E[Y_0] \geq R_0\right\}.
\end{equation*} Hence, for contract $\xi\in \Xi$, the 2BSDEJ representation in Theorem~\ref{thm:agent-2BSDEJ} allows us to
parametrize admissible contracts by $(Y_0,Z,U,K)\in \mathbb Y_0\times \mathbb H^p_0(\mathbb{F}^{\mathcal{P}} , \mathcal{P}) \times \mathbb J^p_0(\mathbb{F}^{\mathcal{P}} , \mathcal{P}) \times \mathbb K^p_0(\mathbb{F}^{\mathcal{P}} , \mathcal{P})$ such that $K$ satisfies \eqref{eq:min.con}. In forward form,
$$
\begin{aligned}
Y_T^{Y_0,Z,U,K}
&=
Y_0
-\int_0^T
G^\star(s,X_s,Y_s,Z_s,U_s;\widehat\sigma_s)ds
+\int_0^T Z_s\cdot dX_s^{c,\mathbb P}
\\
&\quad
+\int_0^T\int_{\mathbb R^n\setminus\{0\}}
U_s(\chi)(\mu_X-\Lambda^0)(ds,d\chi)
+\int_0^T dK_s^{\mathbb P},
\end{aligned}
$$
and hence
$$
    \xi=(U^A)^{-1}\big(Y_T^{Y_0,Z,U,K}-F^A(X_T)\big).
$$
For fixed $Y_0\in \mathbb Y_0$ and let $\mathcal K_{Y_0}$ be the set of processes $(Z,U,K)\in \mathbb H^p_0(\mathbb{F}^{\mathcal{P}} , \mathcal{P}) \times \mathbb J^p_0(\mathbb{F}^{\mathcal{P}} , \mathcal{P}) \times \mathbb K^p_0(\mathbb{F}^{\mathcal{P}} , \mathcal{P})$ such that $(U^A)^{-1}\big(Y_T^{Y_0,Z,U,K}-F^A(X_T)\big)\in \Xi$ and $K$ satsify \eqref{eq:min.con}. The Principal may now optimize over $(Y_0,Z,U,K)\in\mathbb Y_0\times \mathcal K_{Y_0}$, with the participation
constraint $Y_0\ge R_0$.\vspace{0.5em}

Let $\widehat{\mathcal A}$ denote the set of agent best-response selectors.
For readability, we state the PDE below for a fixed measurable selector
$\widehat\alpha\in\widehat{\mathcal A}$. If several selectors are retained,
one adds the outer supremum over $\widehat\alpha\in\widehat{\mathcal A}$ in
the Principal's Hamiltonian.\footnote{For the sake of simplicity, we thus assume that $\hat{\mathcal A}$ is a singleton.} The reduced Principal problem is
$$
\begin{aligned}
V_0^P
:=
\sup_{Y_0\ge R_0}
\sup_{(Z,U,K)\in\mathcal K_{Y_0}}
\inf_{(\mathbb P,\eta)\in\mathcal H^{(\widehat\alpha)}}
\mathbb E^{\mathbb P}\Big[
&F^P(X_T)
-U_A^{-1}\big(Y_T^{Y_0,Z,U,K}-F^A(X_T)\big)
\\
&-\int_0^T C^P(s,X_s,\widehat\alpha_s,\eta_s)ds
\Big].
\end{aligned}
$$

The key obstacle is that the usual dynamic programming principle is not
available in this non-dominated sup-inf setting. We therefore use stochastic
Perron's method. As in \cite{vuca}, the nondecreasing process $K$ is
regularized through an auxiliary matrix process $\Gamma$:
$$
K_s^{z,u,\gamma}
=
\int_t^s
\left(
G^\star(r,X_r,Y_r,z,u;\widehat\sigma_r)
+\frac12\mathrm{Tr}(\widehat\sigma_r\gamma)
-H(r,X_r,Y_r,z,u;\gamma)
\right)dr.
$$
Here $\Gamma$ is an auxiliary regularization variable, not an elementary
control. The elementary Principal controls are $(Z,U,K)$.

For $t\le s\le T$, the controlled Markovian state is
\begin{equation}\label{eq:system}
\begin{cases}
dX_s
\!=\!
b(s,X_{s-};\widehat\alpha_s,\eta_s)ds
+\sigma(s,X_{s-},\eta_s)dW_s^{(\widehat\alpha,\eta)}
\!+\!\displaystyle\int_E \beta(s,X_{s-},e)
\widetilde\mu^{(\widehat\alpha,\eta)}(ds,de),
\\ [1mm]
dY_s
=
\bigg[
\displaystyle\int_E
\big(U_s(\beta(s,X_{s-},e))-Z_s\cdot\beta(s,X_{s-},e)\big)
\bigl(\kappa_s(e;\widehat\alpha_s,\eta_s)-1\bigr)\nu_s^0(de)\\
\,+Z_s\!\cdot\! b(s,X_{s-};\widehat\alpha_s,\eta_s)-G^\star(s,X_{s-},Y_s,Z_s,U_s;\widehat\sigma_s)
\bigg]ds
\!+\!Z_s^\top\sigma(s,X_{s-},\eta_s)dW_s^{(\widehat\alpha,\eta)}\\
\,+dK_s
+\displaystyle\int_E U_s(\beta(s,X_{s-},e))
\widetilde\mu^{(\widehat\alpha,\eta)}(ds,de),
\\ [1mm]
X_t=x,\qquad Y_t=y.
\end{cases}
\end{equation}

\subsection{The Principal's integro-HJBI equation}
\label{subsec:principal-hjbi}
We now introduce the integro-partial HJB-Isaacs PDE associated with the optimization $(2Mm-\sigma)$ reduced to $V_0^P(x,Y_0)$. We refer to \cite{oksendal2019applied,biswas2012zero} for more details about the derivation of this HJB-Isaacs PDE. Let $v:[0,T]\times\mathbb R^n\times\mathbb R\to\mathbb R$ be smooth and set
$$
\Sigma(t,x,\eta):=\sigma\sigma^\top(t,x,\eta),
\qquad
\nu_t^{(\widehat\alpha,\eta)}(de)
:=
\kappa_t(e;\widehat\alpha(t,x,z,u),\eta)\nu_t^0(de).
$$
For controls $(z,u,\gamma)$ and Nature's control $\eta$, define
\begin{align*}
b_Y^{z,u,\gamma,\eta}(t,x,y)
:=
z\cdot b(t,x;\widehat\alpha(t,x,z,u),\eta)
+
\frac12\mathrm{Tr}\big(\Sigma(t,x,\eta)\gamma\big)
-
H(t,x,y,z,u;\gamma)\\
+\displaystyle\int_E
\big(u(\beta(t,x,e))-z\cdot\beta(t,x,e)\big)
\bigl(\kappa_t(e;\widehat\alpha(t,x,z,u),\eta)-1\bigr)\nu_t^0(de)
\end{align*}
$$
b^{z,u,\gamma,\eta}(t,x,y)
:=
\begin{pmatrix}
b(t,x;\widehat\alpha(t,x,z,u),\eta)\\
b_Y^{z,u,\gamma,\eta}(t,x,y)
\end{pmatrix},
$$
and
$$
\mathcal C^{z,\eta}(t,x)
:=
\begin{pmatrix}
\Sigma(t,x,\eta) & \Sigma(t,x,\eta)z\\
z^\top\Sigma(t,x,\eta) & z^\top\Sigma(t,x,\eta)z
\end{pmatrix}.
$$
The positive off-diagonal block is consistent with the convention
$dY^c=z^\top\sigma dW$. The controlled generator is
$$
\begin{aligned}
\mathcal L^{z,u,\gamma,\eta}v(t,x,y)
&=
b^{z,u,\gamma,\eta}(t,x,y)\cdot\nabla v(t,x,y)
+\frac12\mathrm{Tr}\big(\mathcal C^{z,\eta}(t,x)D^2v(t,x,y)\big)
\\
&\quad
+\int_E
\Big[
v\big(t,x+\beta(t,x,e),y+u(\beta(t,x,e))\big)-v(t,x,y)
\\
&\qquad
-D_xv(t,x,y)\cdot\beta(t,x,e)
-\partial_yv(t,x,y)u(\beta(t,x,e))
\Big]\nu_t^{(\widehat\alpha,\eta)}(de).
\end{aligned}
$$
The HJBI operator is
\begin{equation*}\label{eq:HJBI-operator-ref}
\mathcal Q^\star[v](t,x,y)
:=
\sup_{(z,u,\gamma)}
\inf_{\eta\in H}
\left\{
\mathcal L^{z,u,\gamma,\eta}v(t,x,y)
-
C^P(t,x,\widehat\alpha(t,x,z,u),\eta)
\right\}.
\end{equation*}
It leads to the following PDE
\begin{equation}\label{eq:HJBI-ref}
\begin{cases}
-\partial_t v(t,x,y)-\mathcal Q^\star[v](t,x,y)=0,
& (t,x,y)\in[0,T)\times\mathbb R^n\times\mathbb R,\\ [1mm]
v(T,x,y)=F^P(x)-U_A^{-1}\big(y-F^A(x)\big),
& (x,y)\in\mathbb R^n\times\mathbb R.
\end{cases}
\end{equation}
The Principal optimization thus becomes
$$
    V_0^P(x_0)=\sup_{Y_0\ge R_0}v(0,x_0,Y_0).
$$

\subsection{Elementary controls and Perron standing assumptions}
\label{subsec:elem-controls}

Stochastic Perron's method is carried out on elementary strategies. In this
section, $V^P$ denotes the elementary-control Perron value. Under the usual
density/no-loss condition for elementary controls, this coincides with the
unrestricted Principal value. For $s\in[t,T]$, let
$
\mathcal B_s^t:=\sigma\big((X_r,Y_r):t\le r\le s\big)$. 
\begin{definition}[Stopping rule]
A stopping rule is a map
$\tau:D([t,T];\mathbb R^n\times\mathbb R)\to[t,T]$ such that
$\{\tau\le s\}\in\mathcal B_s^t$ for all $s\in[t,T]$.
\end{definition}
\begin{definition}[Elementary controls]
Fix a stopping rule $\tau$.
A Principal elementary control on $[\tau,T]$ is a triple $(Z,U,K)$ such that,
on a finite stopping-time grid $\tau=\tau_0\le\cdots\le\tau_N=T$,
$$
Z_s=\sum_{i=1}^N z_i{\bf 1}_{(\tau_{i-1},\tau_i]}(s),
\qquad
U_s=\sum_{i=1}^N u_i{\bf 1}_{(\tau_{i-1},\tau_i]}(s),
$$
and
$$
K_s=\int_\tau^s k_rdr,
\qquad
k_r=\sum_{i=1}^N k_i{\bf 1}_{(\tau_{i-1},\tau_i]}(r),
\qquad k_i\ge0.
$$
Here $(z_i,u_i,k_i)$ are measurable with respect to
$\mathcal B_{\tau_{i-1}}^t$. We denote the set of such controls by
$\mathfrak K(t,\tau)$. Nature's elementary controls are defined analogously by
$$
\eta_s=\sum_{i=1}^N h_i{\bf 1}_{(\tau_{i-1},\tau_i]}(s),
\qquad h_i\in H,
$$
and the corresponding set is denoted by $\mathfrak H(t,\tau).$
\end{definition}

We define the elementary control set as \begin{equation*}
    \mathfrak K_{Y_0} = \mathcal K_{Y_0}\cap \mathfrak K
\end{equation*}

Given a best-response selector \(\widehat\alpha(t,x,z,u)\), we write
\(\mathfrak P^{(\widehat\alpha)}(t,\tau)\) for the collection of admissible
continuation prior-control pairs after the stopping rule \(\tau\):
$$
\mathfrak P^{(\widehat\alpha)}(t,\tau)
:=
\Big\{
(\mathbb P,\eta):\
\eta\in\mathfrak H(t,\tau),\ 
\mathbb P\in\mathcal P^{(\widehat\alpha)}(t,\tau)
\text{ is consistent with }\eta
\Big\}.
$$
When \(\tau=t\), we write \(\mathfrak P^{(\widehat\alpha)}(t,t)\) for the
full-path admissible pairs on \([t,T]\).

\begin{assumption}[Standing assumptions for the Perron argument]
\label{ass:perron-standing}
The following conditions hold.

\begin{enumerate}
\item[(i)] \emph{Local Hamiltonian localization.}
For every $\phi\in C^{1,2}$ and every $(t,x,y)$, the supremum in
$\mathcal Q^\star[\phi](t,x,y)$ may be restricted to a localized compact set
$\mathcal K_R(t,x,y)\subset\mathbb R^n\times\mathcal L_\nu^p\times\mathbb S^n$:
$$
\mathcal Q^\star[\phi](\cdot)
=
\sup_{(z,u,\gamma)\in\mathcal K_R(\cdot)}
\inf_{\eta\in H}
\left\{
\mathcal L^{z,u,\gamma,\eta}\phi(\cdot)
-
C^P(t,x,\widehat\alpha(t,x,z,u),\eta)
\right\}.
$$

\item[(ii)] \emph{Finite-activity jumps and moments.}
For every localized set $\mathcal K_R$, there exists $C_R>0$ such that
$$
\sup_{(z,u,\gamma)\in\mathcal K_R,\ \eta\in H}
\nu_t^{(\widehat\alpha,\eta)}(E)
\le C_R,
\qquad t\in[0,T].
$$
Moreover, for some integer $q\ge1$ dominating the polynomial growth of the
terminal payoff and semisolutions, we have
$$
\sup_{\substack{(z,u,\gamma)\in\mathcal K_R\\\eta\in H}}
\int_E
\left(
|\beta(t,x,e)|^{2q}
+
|u(\beta(t,x,e))|^{2q}
\right)
\nu_t^{(\widehat\alpha,\eta)}(de)
\!\le\!
C_R(1+|x|^{2q}+|y|^{2q}).
$$

\item[(iii)] \emph{Stability under elementary conditioning and pasting.}
Let $\theta$ be a stopping rule and let
$(Z^0,U^0,K^0)\in\mathfrak K_{Y_0}(t,t)$, $\eta^0\in\mathfrak H(t,t)$
be elementary controls used up to $\theta$. Let
$(A_j)_{j=1}^N\subset\mathcal B_\theta^t$ be a finite partition and, for each
$j$, let
$$
    (Z^j,U^j,K^j)\in\mathfrak K_{y_0}(t,\theta),
    \qquad
    \eta^j\in\mathfrak H(t,\theta)
$$
be elementary continuation controls after $\theta$. Then the pasted controls
defined by
$$
Z_s
=
Z_s^0{\bf 1}_{\{s\le\theta\}}
+
\sum_{j=1}^N{\bf 1}_{A_j}Z_s^j{\bf 1}_{\{s>\theta\}},
$$
$$
U_s
=
U_s^0{\bf 1}_{\{s\le\theta\}}
+
\sum_{j=1}^N{\bf 1}_{A_j}U_s^j{\bf 1}_{\{s>\theta\}},
$$
and
$$
k_s
=
k_s^0{\bf 1}_{\{s\le\theta\}}
+
\sum_{j=1}^N{\bf 1}_{A_j}k_s^j{\bf 1}_{\{s>\theta\}},
\qquad
K_s:=\int_t^s k_rdr,
$$
belong to $\mathfrak K_{y_0}(t,t)$. Similarly, the pasted Nature control
$$
\eta_s
=
\eta_s^0{\bf 1}_{\{s\le\theta\}}
+
\sum_{j=1}^N{\bf 1}_{A_j}\eta_s^j{\bf 1}_{\{s>\theta\}}
$$
belongs to $\mathfrak H(t,t)$. Moreover, the induced family of laws is stable under the same pasting operation:
if $(\mathbb P^0,\eta^0)\in\mathfrak P^{(\widehat\alpha)}(t,t)$ and, on each
$A_j$, the continuation law is generated by an admissible pair
$(\mathbb P^j,\eta^j)\in\mathfrak P^{(\widehat\alpha)}(t,\theta)$, then the
law obtained by following $\mathbb P^0$ up to $\theta$ and $\mathbb P^j$
after $\theta$ on $A_j$ again belongs to
$\mathfrak P^{(\widehat\alpha)}(t,t)$. Equivalently, the corresponding pasted
Girsanov kernels $(\zeta,\kappa)$ remain admissible and the associated
Doléans--Dade exponential is a true martingale.

\end{enumerate}
\end{assumption}

\subsection{Stochastic semisolutions and Perron envelopes}
\label{subsec:perron}
We now turn to the definitions of stochastic super and sub solution to \eqref{eq:HJBI-ref}. Let
$$
g(x,y):=F^P(x)-U_A^{-1}\big(y-F^A(x)\big).
$$

\begin{definition}[Stochastic semisolutions of \eqref{eq:HJBI-ref}]
\label{def:semisolution}
Let $Y_0\in \mathbb Y_0$. Let $
v \colon [0,T]\times \mathbb{R}^n \times \mathbb{R} \longrightarrow \mathbb{R}$

\begin{itemize}
\item \textbf{Sub-solution.} $v$ is called a \emph{stochastic sub-solution} of the HJBI equation \eqref{eq:HJBI-ref} if \begin{itemize}
\item[(i-)] $v$ is continuous and 
$$
v(T, x,y)\le g(x,y)
\quad\text{for any }( x,y)\in \mathbb{R}^n\times \mathbb{R},
$$ 
\item[(ii-)] for any $t\in[0,T]$ and for any stopping rule $\tau\in \mathcal{B}^t$, 
there exists an \emph{elementary control} $(\widetilde{Z},\widetilde{U},\widetilde{K})\in \mathfrak K_{Y_0}(t,\tau)$ 
such that for \emph{any} $(Z,U,K)\in \mathfrak K (t,t)$, 
for any $(\mathbb{P},\eta)\in \mathfrak{P}^{\hat \alpha}(t,t)$ 
and every stopping rule $\rho \in \mathcal{B}^t$ with $\tau \le \rho \le T$ we have
$$
v\bigl(\tau', X^{(\tau)}_{\tau'},Y^{(\tau)}_{\tau'}\bigr)\le 
\mathbb{E}^{\mathbb{P}}\Bigl[
v\bigl(\rho', X^{(\tau)}_{\rho'},Y^{(\tau)}_{\rho'}\bigr)-
\int_{\tau'}^{\rho'}
C^P(s,X_s,\widehat\alpha_s,\eta_s)ds\Big|
\mathcal{F}^{t}_{\tau'}\Bigr]
$$
where, for any $( x,y,\omega)\in\mathbb{R}^n\times \Omega$,
$$
 X^{(\tau)} := 
 X^{t, x,(Z,U,K)\otimes_{\tau}(\widetilde{Z},\widetilde{U},\widetilde{K}),\eta},
\quad
Y^{(\tau)} :=
Y^{t,y,(Z,U,K)\otimes_{\tau}(\widetilde{Z},\widetilde{U},\widetilde{K}),\eta},
$$ where $ X^{t, x,(Z,U,K)\otimes_{\tau}(\widetilde{Z},\widetilde{U},\widetilde{K}),\eta},Y^{t,y,(Z,U,K)\otimes_{\tau}(\widetilde{Z},\widetilde{U},\widetilde{K}),\eta} $ denotes the solution to the controlled system \eqref{eq:system}, with concatenated elementary strategies control $(\widetilde{Z},\widetilde{U},\widetilde{K})$ starting with $(Z,U,K)$ at time $t$, see \cite[Definition 3.1]{perronsirbu}
$$
\tau'(\omega)
:=
\tau\bigl(
 X^{t, x,(Z,U,K)\otimes_{\tau}(\widetilde{Z},\widetilde{U},\widetilde{K}),\eta}(\omega),
Y^{t,y,(Z,U,K)\otimes_{\tau}(\widetilde{Z},\widetilde{U},\widetilde{K}),\eta}(\omega)
\bigr),
$$
$$
\rho'(\omega)
:=
\rho\bigl(
 X^{t, x,(Z,U,K)\otimes_{\tau}(\widetilde{Z},\widetilde{U},\widetilde{K}),\eta}(\omega),
Y^{t,y,(Z,U,K)\otimes_{\tau}(\widetilde{Z},\widetilde{U},\widetilde{K}),\eta}(\omega)
\bigr).
$$

\end{itemize} We denote by $\mathcal{V}^{-}$ the set of all such \emph{stochastic sub-solutions} to \eqref{eq:HJBI-ref} .

\vspace{1em}

\item \textbf{Super-solution.} $v$ is a \emph{stochastic super-solution} of the HJBI equation \eqref{eq:HJBI-ref}  if
\begin{itemize}
\item[(i+)] $v$ is continuous and 
$$
v(T, x,y)\ge F^P( x)-U_A^{-1}(y-F^A( x))
\quad\text{for any }( x,y)\in \mathbb{R}^n\times \mathbb{R},
$$
\item[(ii+)] for any $t\in [0,T]$, for any stopping rule $\tau\in \mathcal{B}^t$ 
and for any $(Z,U,K)\in \mathfrak K_{Y_0}(t,\tau)$, 
there exists an \emph{elementary control} $(\widetilde{\mathbb{P}},\widetilde{\eta})\in \mathfrak P^{\hat\alpha}(t,\tau)$ 
such that for every $\eta \in \mathfrak H(t,t)$ satisfying $(\widetilde{\mathbb{P}},\eta)\in \mathfrak P^{\hat\alpha}(t,t)$ 
and for every stopping rule $\rho\in \mathcal{B}^t$ with $\tau \le\rho \le T$, we have
$$
v\bigl(\tau', X^{(\tau)}_{\tau'},Y^{(\tau)}_{\tau'}\bigr)\ge 
\mathbb{E}^{\widehat{\mathbb{P}}}\Bigl[
v\bigl(\rho', X^{(\tau)}_{\rho'},Y^{(\tau)}_{\rho'}\bigr)-\int_{\tau'}^{\rho'}
C^P(s,X_s,\widehat\alpha_s,\eta_s)ds\Big|
\mathcal{F}^{t}_{\tau'}\Bigr]
$$
where, for any $( x,y,\omega)\in\mathbb{R}^n\times \Omega$,
$$
 X^{(\tau)} := 
 X^{t, x,Z,U,K,\eta\otimes_\tau\tilde\eta},
\quad
Y^{(\tau)} :=
Y^{t, x,Z,U,K,\eta\otimes_\tau\tilde\eta},
$$ 
$$
\tau'(\omega)
:=
\tau\bigl(
 X^{t, x,Z,U,K,\eta\otimes_\tau\tilde\eta}(\omega),
Y^{t, x,Z,U,K,\eta\otimes_\tau\tilde\eta}(\omega)
\bigr),
$$
$$
\rho'(\omega)
:=
\rho\bigl(
 X^{t, x,Z,U,K,\eta\otimes_\tau\tilde\eta}(\omega),
Y^{t, x,Z,U,K,\eta\otimes_\tau\tilde\eta}(\omega)
\bigr).
$$

\end{itemize}
We denote by $\mathcal{V}^{+}$ the set of all such \emph{stochastic super-solutions} to \eqref{eq:HJBI-ref} .
\end{itemize}
\end{definition}

\begin{assumption}The sets $\mathcal{V}^{+}$ and $\mathcal{V}^{-}$ are non-empty.\end{assumption}

As explained in \cite{perrondynkin,bayraktar2012stochastic}, the set $\mathcal{V}^{+}$ is trivially non-empty if $U_{P}$ is bounded above, whereas $\mathcal{V}^{-}$ is non-empty if $U_{P}$ is bounded below. We now follow Perron's method as in \cite{vuca}. Define \begin{equation}
\label{eq:perron-bracketing}
    v^{-}:=\sup_{v \in \mathcal{V}^{-}} v,\quad v^{+}:=\inf_{v \in \mathcal{V}^{+}} v.
\end{equation}

\begin{theorem}[Perron viscosity envelopes]
\label{thm:viscosity}
Under Assumption~\ref{ass:perron-standing}, $v^{-}$ is a lower semicontinuous
viscosity supersolution of \eqref{eq:HJBI-ref} on
$[0,T)\times\mathbb R^n\times\mathbb R$, and $v^{+}$ is an upper
semicontinuous viscosity subsolution on the same domain.
\end{theorem}

\begin{remark}
The proof follows \cite[Thm.~3.1]{BayraktarLi2016}, \cite[Thm.~3.5]{perronsirbu} and \cite[Thm.~4.1]{vuca}, with the only structural modification on the global coercive patch needed to handle jump overshoots, together with the inclusion of the jump integrand $U$ in the elementary Principal control. Full details are provided in Appendix \ref{app:viscosity_proof}.
\end{remark}

\subsection{Comparison and identification}
\label{subsec:comparison-growth}

While Theorem~\ref{thm:viscosity} together with 
 \eqref{eq:perron-bracketing} enables us to derive the following inequalities 
$$
    \bar v^{-}\le V^P\le \bar v^{+},
$$ we need a comparison principle in order to fully characterize the value function as a viscosity solution to the PDE \eqref{eq:HJBI-ref}. For $m\ge1$, let $\mathrm{USC}_m$ and $\mathrm{LSC}_m$ denote the upper
and lower semicontinuous functions on
$[0,T]\times\mathbb R^n\times\mathbb R$ satisfying
$$
    |w(t,x,y)|\le C(1+|x|^m+|y|^m)
$$
for some $C>0$.

\begin{assumption}[Comparison in a polynomial-growth class]
\label{ass:comparison-principle}
There exists $m\ge1$ such that comparison holds for \eqref{eq:HJBI-ref}:
if $u\in\mathrm{USC}_m$ is a viscosity subsolution, $w\in\mathrm{LSC}_m$ is
a viscosity supersolution, and $u(T,\cdot,\cdot)\le w(T,\cdot,\cdot)$, then
$u\le w$ on $[0,T]\times\mathbb R^n\times\mathbb R$.
\end{assumption}

\begin{remark}
Note that sufficient conditions ensuring that the comparison principle above is satisfied include polynomial growth of $g$, continuity of
the selected feedback $\widehat\alpha$, local Lipschitz and at most linear
growth of $b,\mathcal C$, properness and degenerate ellipticity of
the operator $\mathcal Q^*$, localization of $(z,u,\gamma)$, and (ii) in Assumption~\ref{ass:perron-standing}. These conditions place
\eqref{eq:HJBI-ref} within the usual comparison frameworks for fully nonlinear
integro-differential Bellman--Isaacs equations; see
\cite{BarlesImbert2008,JakobsenKarlsen2005,JakobsenKarlsen2006}. In the cyber
application, the epidemic coordinates are bounded proportions, controls belong to compact sets, jump-size laws have finite support, jump intensities are bounded on the computational domain, and terminal payoffs have polynomial growth.
\end{remark}

\begin{corollary}[Identification under comparison]
\label{cor:identification-comparison}
Under Assumption~\ref{ass:comparison-principle}. We define
$$
\bar v^{-}(t,x,y)
:=
\begin{cases}
v^{-}(t,x,y),&t<T,\\
g(x,y),&t=T,
\end{cases}
\qquad
\bar v^{+}(t,x,y)
:=
\begin{cases}
v^{+}(t,x,y),&t<T,\\
g(x,y),&t=T.
\end{cases}
$$
If $\bar v^{-}\in\mathrm{LSC}_m$ and
$\bar v^{+}\in\mathrm{USC}_m$, then
$$
    \bar v^{-}=\bar v^{+}=:v.
$$
Moreover,
$$
    V_t^P(x,y)=v(t,x,y),
    \qquad (t,x,y)\in[0,T]\times\mathbb R^n\times\mathbb R,
$$
and
$$
    V_0^P(x_0)=\sup_{Y_0\ge R_0}v(0,x_0,Y_0).
$$
Thus, the Principal's value is identified with the unique
viscosity solution of \eqref{eq:HJBI-ref} in the chosen polynomial-growth class.
\end{corollary}

\begin{proof}
The relaxed terminal inequalities in Theorem~\ref{thm:viscosity} ensure that
the terminal extensions $\bar v^{-}$ and $\bar v^{+}$ are respectively lower
and upper semicontinuous at $T$. Since $\bar v^{+}(T,\cdot,\cdot)
=\bar v^{-}(T,\cdot,\cdot)=g$, comparison gives
$ \bar v^{+}\le \bar v^{-}.$
On the other hand, the Perron bracketing \eqref{eq:perron-bracketing} gives
$$
    \bar v^{-}\le V^P\le \bar v^{+}.
$$
Hence $\bar v^{-}=V^P=\bar v^{+}$. Uniqueness follows by applying comparison
in both directions to any two polynomial-growth viscosity solutions.
\end{proof}

}

\section{Application in cyber risk management}
\label{sec:cyber}

\subsection{Cyber risk modeling: controlled SIR–price system}

We now turn to the particular cyber risk model we are considering by specifying the dynamic of $X$ with a controlled SIR model and the subsidiary's portfolio evolution. 
We model the computers or electronic devices in the cluster by SIR model, following the construction in \cite{hillairet2024optimal}: 

\hspace{0.15\textwidth}
\begin{minipage}{0.7\textwidth}
    \centering
    \textbf{SIR model for cyber contagion and attacks}
    \vspace{0.1cm}
    \begin{tikzpicture}
    \node [rectangle, draw, rounded corners=4pt, fill=blue!20, text centered] (S) at (0,0) {\bf Susceptible};
    \node [rectangle, draw, rounded corners=4pt, fill=red!20, text centered, dotted, thick] at (4,0) (I) {\bf Infected};
    \node [rectangle, draw, rounded corners=4pt, fill=green!20, text centered, below of=I] (D) at (2,-1.) {\bf Recovered};
    [node distance=1cm] 
    \draw[->, draw=black, thick]  (S) edge[bend left]  node[midway, above] {\bf contagion }(I);
    \draw[->, draw=black, thick] (S) edge[bend right]  node[midway, below] {\bf hacking}(I);
    \draw[->, draw=black, thick] (I) edge[bend left]  node[midway, right] {\bf replacement }(D);
    \draw[->, draw=black, thick]  (S) edge[bend right]  node[midway, left] {\bf protection }(D);
    \end{tikzpicture}

\end{minipage}%

\begin{itemize}
    \parskip = 0.1pt
    \item \textbf{Susceptible (S)}: $ S_t $ denotes the proportion of computers at time $ t $ that are insufficiently protected and not yet infected, making them susceptible to attacks.
    \item \textbf{Infected (I)}: $ I_t $ represents the proportion of infected and corrupted computers at time $ t $ that can potentially contaminate other devices through cyber contagion and interconnectedness.
    \item \textbf{Recovered (R)}: $ R_t $ indicates the proportion of computers at time $ t $ that have either recovered from infection or are protected by antivirus software, rendering them immune to future infections. 
\end{itemize}

Under any admissible $(\alpha,\eta)$, the controlled SIR system is
\begin{equation}\label{eq:SIR}
\begin{cases}
dS_t = \big(-\bar\beta S_t I_t - \alpha_t S_t - \eta_t S_t\big)dt - \widetilde\sigma(t,\eta_t)S_t I_td\widetilde W_t,\\
dI_t = \big(\bar\beta S_t I_t - \rho I_t + \eta_t S_t\big)dt + \widetilde\sigma(t,\eta_t)S_t I_td\widetilde W_t,\\
dR_t=(\rho I_t+\alpha_tS_t)dt,\\
S_t+I_t+R_t=1
\end{cases}
\end{equation}

\paragraph{Transmission and controls.}
The constant $\bar\beta>0$ is the baseline transmission rate and $\rho>0$ is the recovery rate. The hacker’s control $\eta_t\in H$ modulates both the epidemic and volatility, and also affects the portfolio’s volatility and jump intensities. The subsidiary’s (agent’s) control $\alpha_t\in A$ is a protection effort acting on $S$.

\subsection{L–hop modeling of jump sources}

{ To illustrate the L-hop propagation mechanism, we specify two marked Poisson
random measures $N^e(dt,dc)$ and $N^i(dt,dc)$ on $(0,\infty)\times(0,1)$.
The marks $c\in(0,1)$ represent relative cyber-loss severities. Their predictable
compensators are
$$
    \lambda^e(\eta_t)\nu^e(dc)dt,
    \qquad
    \lambda^i(I_t,\alpha_t)\nu^i(dc)dt,
$$
where $\nu^e$ and $\nu^i$ are probability laws on $(0,1)$ describing,
respectively, external and internal jump-size distributions. The price dynamics
are then
\begin{equation}\label{eq:P-Lhop-random}
\frac{dP_{t-}}{P_{t-}}
=
\mu(t,I_t)dt+\sigma_P(t,I_t,\eta_t)dW_t
-\int_{(0,1)} c\,N^e(dt,dc)
-\int_{(0,1)} c\,N^i(dt,dc).
\end{equation}
Thus each external or internal cyber incident instantaneously scales the
portfolio value by a random factor $1-c$. In the final numerical calibration,
$\nu^e$ and $\nu^i$ are chosen as finite three-point distributions; see
Appendix~\ref{appendix:parameter_final}.}

\subsection{Admissible control}

The subsidiary’s effort $\alpha\in\mathfrak A$ is $\mathbb F$–progressively measurable with values in a compact set $A$, and acts through the drift of $S$. For clarity, we decompose the drift of $\mathbf X=(P,S,I)$ as
$$
\mathbf b^\eta(t,\mathbf X_{t-},\eta_t)
:=\begin{pmatrix}
\mu(t,I_t)P_{t-}\\
-\bar\beta S_t I_t-\eta_t S_t\\
\bar\beta S_t I_t+\eta_t S_t-\rho I_t
\end{pmatrix},\qquad
\mathbf b^\alpha(\mathbf X_t;\alpha_t)
:=\begin{pmatrix}
0\\
-\alpha_t S_t\\
0
\end{pmatrix}.
$$
The continuous volatility matrix (two Brownian directions) is
$$
\boldsymbol\sigma(t,\mathbf X_{t-},\eta_t)
=\begin{pmatrix}
\sigma_P(t,I_t,\eta_t)P_{t-} & 0\\
0 & -\widetilde\sigma(t,\eta_t)S_t I_t\\
0 & \widetilde\sigma(t,\eta_t)S_t I_t
\end{pmatrix},
$$ where we take $\ell =2$ for the Brownian motion dimension, with the first component driring the price, the second driving the SIR system.
Combining with the jump part from \eqref{eq:P-Lhop-random}, one can write compactly
\begin{equation*}
    \begin{split}
        d\mathbf X_{t-}
=\Big(\mathbf b^\eta(t,\mathbf X_{t-},\eta_t)+\mathbf b^\alpha(\mathbf X_{t-};\alpha_t)\Big)dt
+\boldsymbol\sigma(t,\mathbf X_{t-},\eta_t)d\mathbf W_t
\\
+\begin{pmatrix}
-P_{t-}\cdot\int_{(0,1)} c\,N^e(dt,dc)
-P_{t-}\cdot\int_{(0,1)} c\,N^i(dt,dc)\\
0\\
0
\end{pmatrix}.
    \end{split}
\end{equation*}

\begin{remark}
    All assumptions needed for boundedness and Lipschitz of coefficients, compact controls, and bounded intensities are satisfied, see Appendix \ref{appendix:parameter_final}, which ensures the controls we consider are admissible.
\end{remark}

{
In this case, the optimal contract for a risk-neutral agent is given by 
$$
    \xi=Y_T^{Y_0,Z,U,K}-F^A(X_T),
$$
where

$$
\begin{aligned}
Y_T^{Y_0,Z,U,K}
&=
Y_0
-\int_0^T
G^\star(s,X_s,Y_s,Z_s,U_s;\widehat\sigma_s)ds
+\int_0^T Z_s\cdot dX_s^{c,\mathbb P}
\\
&\quad
+\int_0^T\int_{\mathbb R^n\setminus\{0\}}
U_s(\chi)(\mu_X-\Lambda^0)(ds,d\chi)
+\int_0^T dK^{Z,U,\Gamma}_s.
\end{aligned}
$$

 The optimal contract is dynamically implemented through the continuation-utility process and can be interpreted through its loading coefficients. The $Y_0$ is calibrated to match the reservation utility of the agent. The $G^\star$ is the certain equivalent utility of the agent seen as the surplus of weath the agent is making by optimizing his own utility. The $z$-component is a continuous payment (a penalty in our case) with respect to the evolution of the canonical process $X=(P,S,I)$: it determines how the agent’s continuation payoff reacts to marginal fluctuations in the observable state variables, such as portfolio performance and the cyber state. The jump term $U$ captures discrete contract adjustments following cyber incidents. The $\Gamma$-process captures exposure to quadratic variation and may be interpreted as variance sharing under ambiguity. In the cyber application, this means the contract can be read as a state-contingent compensation rule indexed by observable proxies $(P,S,I)$, allowing us to study whether compensation loads more strongly on performance preservation or cyber resilience in highly infected regimes.
}

\subsection{Numerical simulation and results}
\label{subsec:numerics-setup}
\begin{figure}[H]
     \centering
     \begin{subfigure}[b]{0.47\textwidth}
         \centering
         \includegraphics[width=\textwidth]{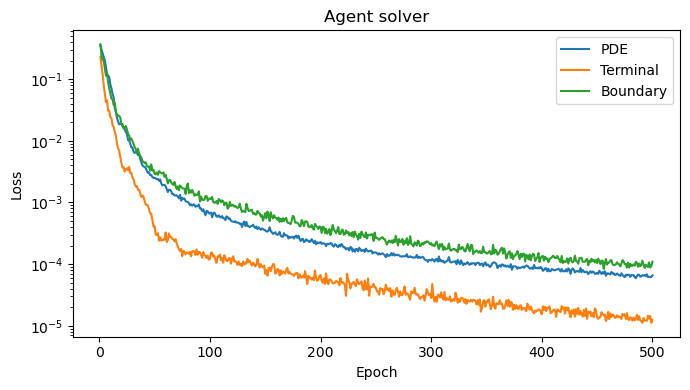}
         \caption{Agent's value without contract}
         \label{fig:agent_conv}
     \end{subfigure}\\
     \hfill
     \begin{subfigure}[b]{0.47\textwidth}
         \centering
         \includegraphics[width=\textwidth]{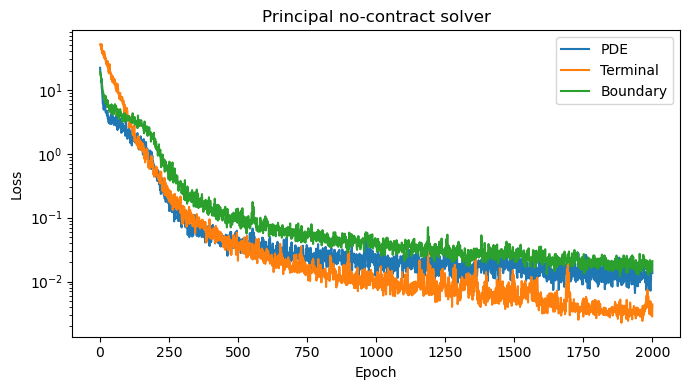}
         \caption{Principal's value without contract}
         \label{fig:principal_nc_conv}
     \end{subfigure}
     \hfill
     \begin{subfigure}[b]{0.47\textwidth}
         \centering
         \includegraphics[width=\textwidth]{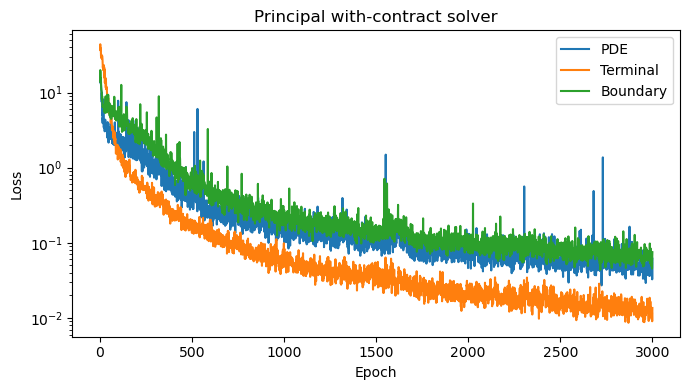}
         \caption{Principal's value with contract}
         \label{fig:principal_wc_conv}
     \end{subfigure}
        \caption{Convergence Graphs associated with Algorithm E.1 and E.2}
        \label{fig:Conv_Graph}
\end{figure}

\begin{figure}
    \centering
     \begin{subfigure}[b]{0.43\textwidth}
         \begin{subfigure}[b]{\textwidth}
             \centering
             \includegraphics[width=\textwidth]{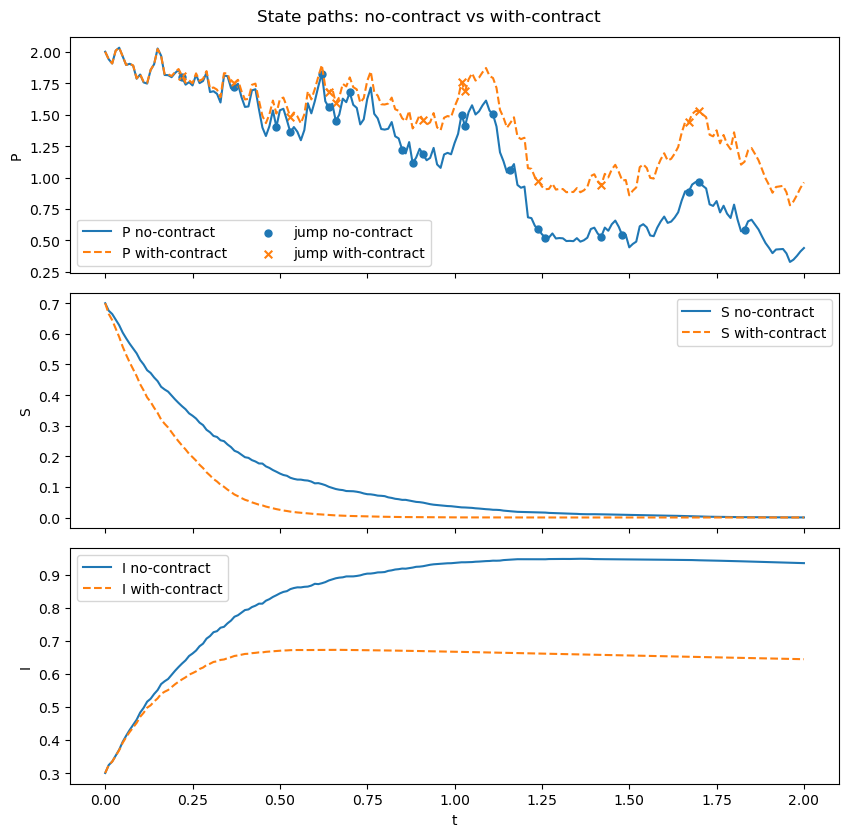}\hfill
            \caption{one example P,S,I path}
            \label{fig:one_path}
        \end{subfigure}
        \begin{subfigure}[b]{\textwidth}
             \centering
             \includegraphics[width=\textwidth]{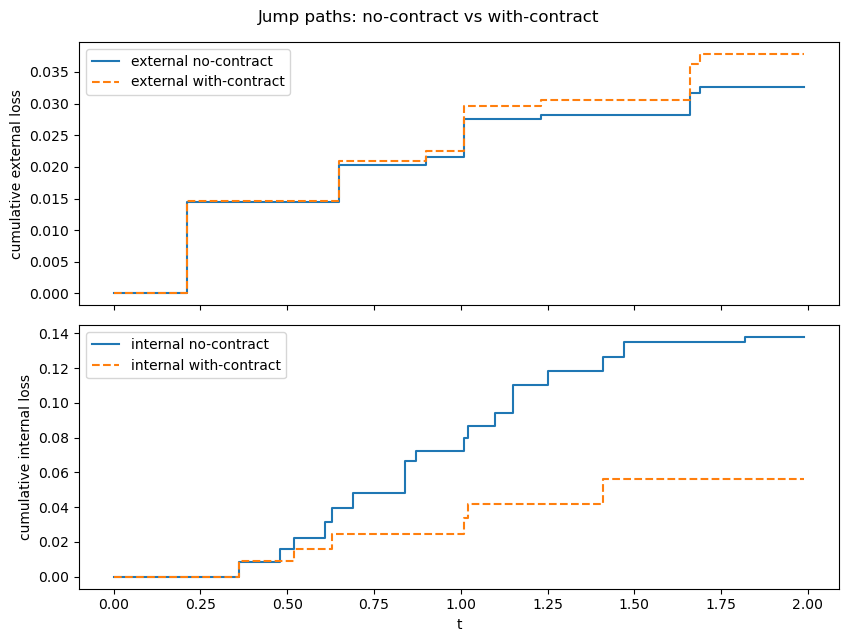}\hfill
            \caption{Jump path}
            \label{fig:one_jump}
        \end{subfigure}
     \end{subfigure}
     \hfill
     \begin{subfigure}[b]{0.54\textwidth}
        \begin{subfigure}[b]{\textwidth}
             \centering
            \includegraphics[width=1.0\linewidth]{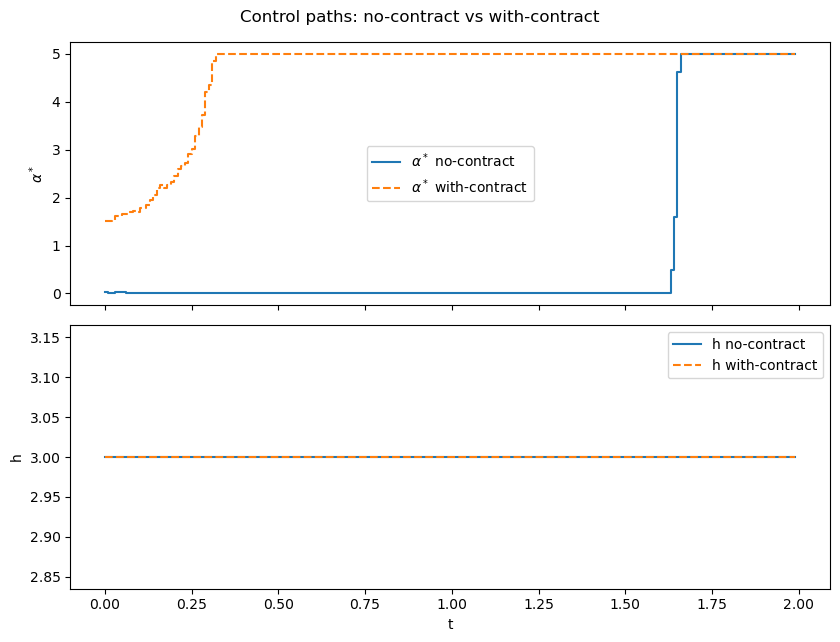}
            \caption{One path control}
            \label{fig:one_control}
        \end{subfigure}\\
        \begin{subfigure}[b]{\textwidth}
             \centering
            \includegraphics[width=\linewidth]{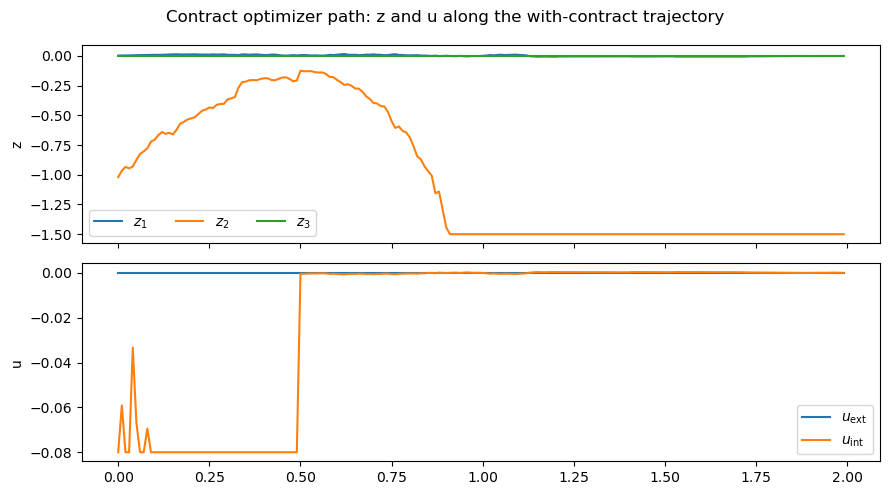}
            \caption{One path optimizer}
            \label{fig:one_optimizer}
        \end{subfigure}
     \end{subfigure}
     \caption{Single-path simulation results}
    \label{fig:One_Path_Graphs}
\end{figure}

\begin{figure}
     \begin{subfigure}[b]{0.43\textwidth}
         \begin{subfigure}[b]{\textwidth}
             \centering
             \includegraphics[width=\textwidth]{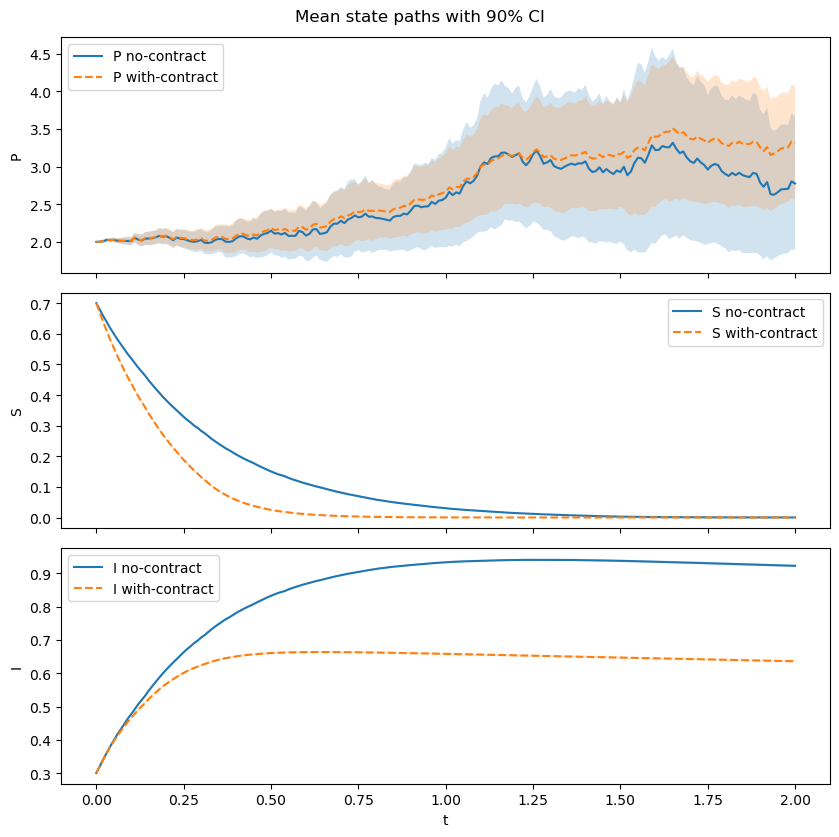}
             \caption{32 paths of P,S,I with 90\% CI}
             \label{fig:32_paths}
         \end{subfigure}\\
         \begin{subfigure}[b]{\textwidth}
             \centering
             \includegraphics[width=\textwidth]{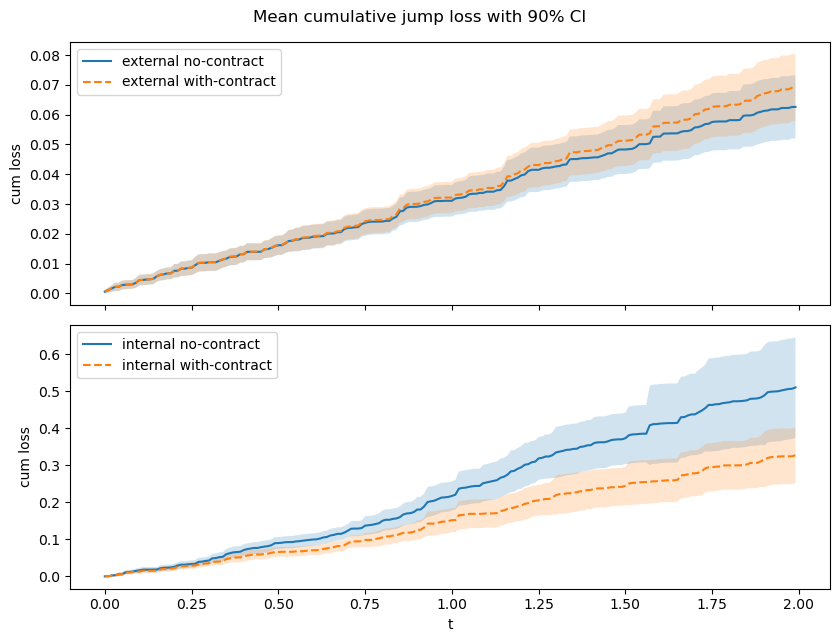}
             \caption{32 jump paths}
             \label{fig:32_jump}
         \end{subfigure}
     \end{subfigure}
     \hfill
     \begin{subfigure}[b]{0.54\textwidth}
        \begin{subfigure}[b]{\textwidth}
             \centering
            \includegraphics[width=\linewidth]{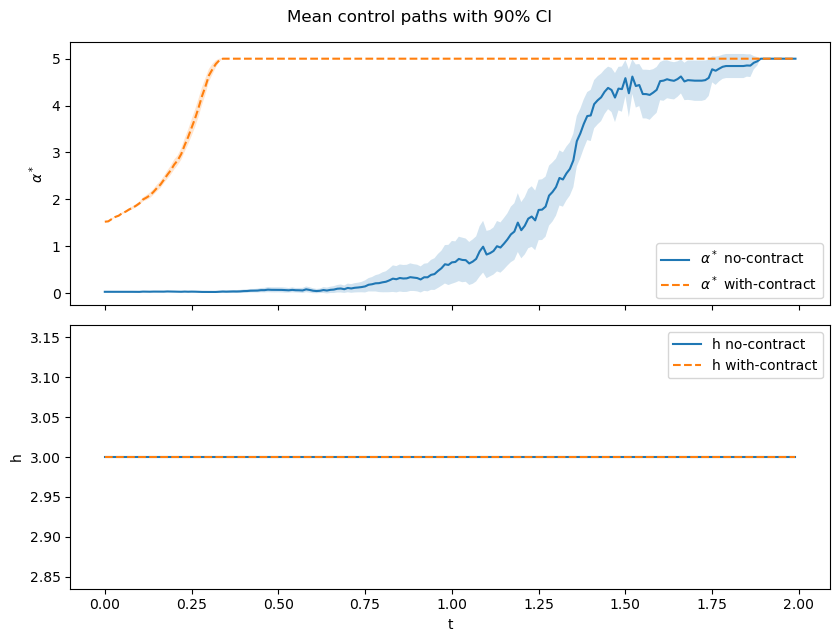}
            \caption{32 paths control}
            \label{fig:32_control}
        \end{subfigure}\\
        \begin{subfigure}[b]{\textwidth}
             \centering
            \includegraphics[width=\linewidth]{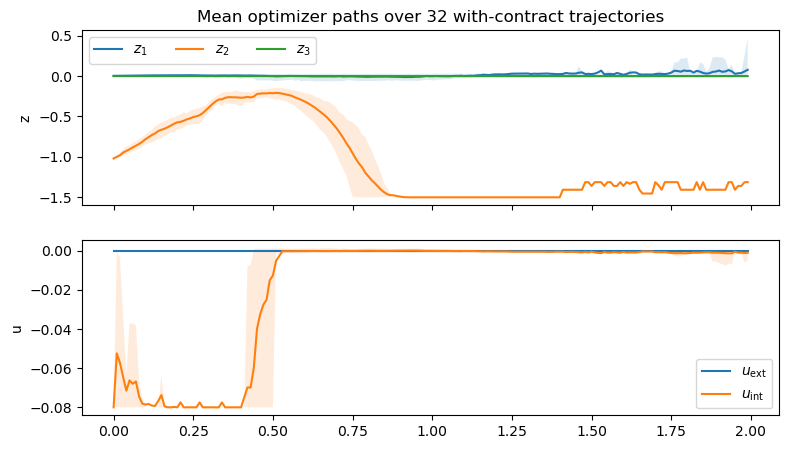}
            \caption{32 paths optimizer}
            \label{fig:32_optimizer}
        \end{subfigure}
     \end{subfigure}
        \caption{Monte Carlo 32 paths simulations}
        \label{fig:32_Path_Graphs}
\end{figure}

\begin{figure}
    \centering
     \begin{subfigure}[b]{0.47\textwidth}
         \centering
         \includegraphics[width=\textwidth]{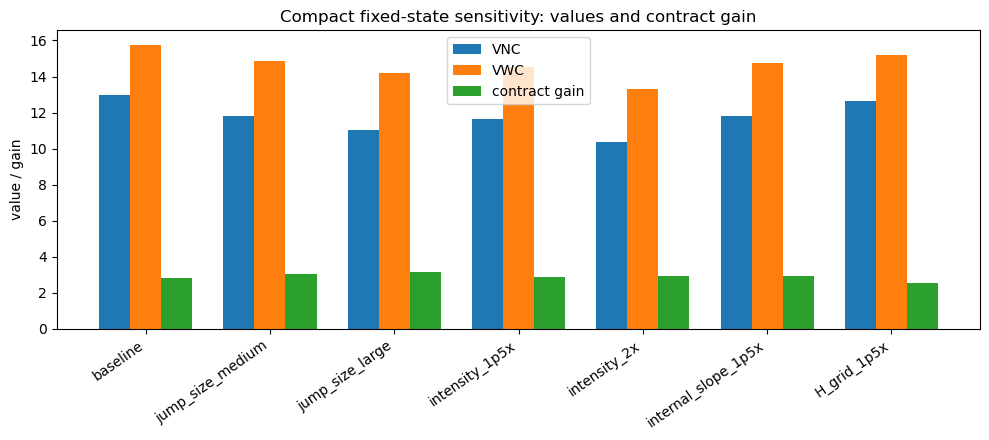}\hfill
        \caption{Principal values and contract gains under different scenarios}
     \end{subfigure}
     \hfill
     \begin{subfigure}[b]{0.47\textwidth}
         \centering
         \includegraphics[width=\textwidth]{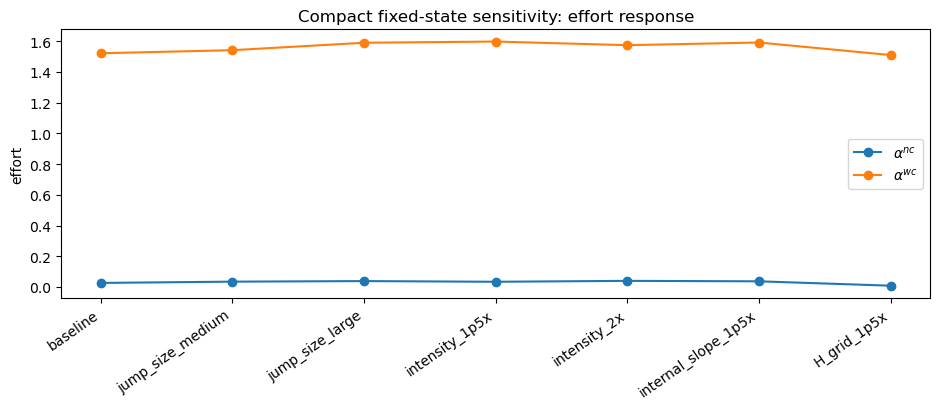}\hfill
        \caption{Average agent's starting effort under different scenarios}
     \end{subfigure}
        \caption{Sensitivity analysis}
        \label{fig:Sen_Ana}
\end{figure}

\begin{figure}
    \centering
    \begin{subfigure}[b]{1\textwidth}
        \centering
        \includegraphics[width=1\textwidth]{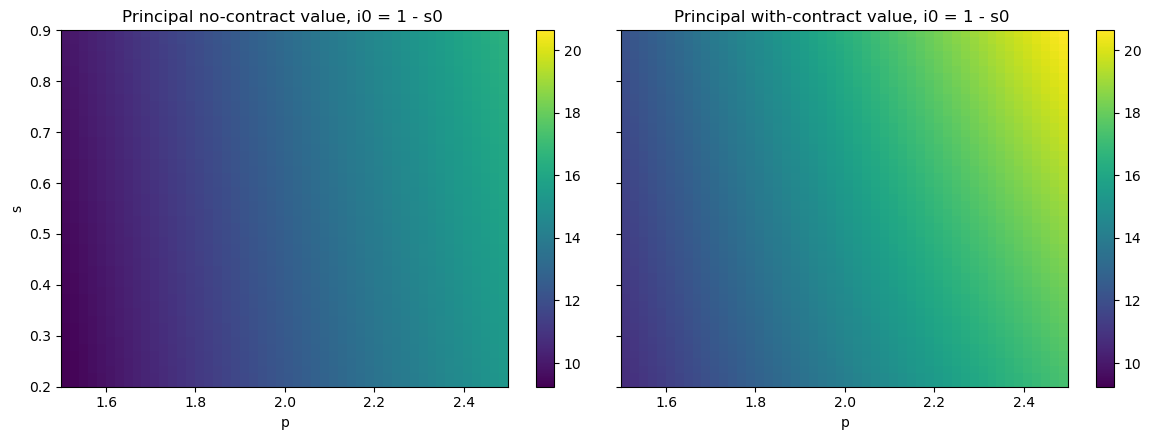}\hfill
        \caption{Principal values without and with contract under different initializations of $s_0$ and $p_0$}
    \end{subfigure}\\
    \begin{subfigure}[b]{0.47\textwidth}
         \centering
         \includegraphics[width=\textwidth]{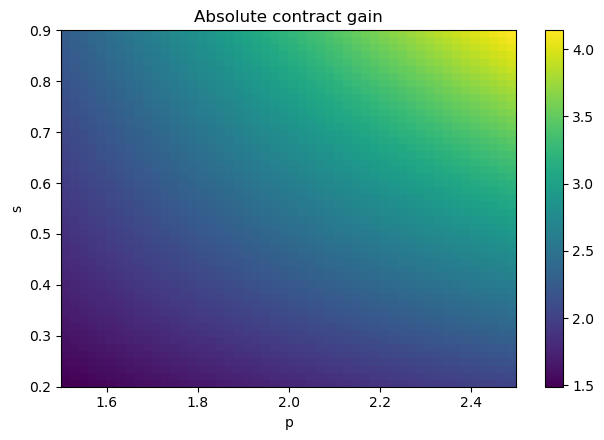}\hfill
        \caption{Principal values obtained from the contract under different initializations of $s_0$ and $p_0$}
     \end{subfigure}
     \caption{Principal values analysis}
     \label{fig:principal_values}
\end{figure}

{
We solve the HJB/HJBI equations with physics-informed neural networks (PINNs) in the spirit of DGM \cite{sirignano2018dgm}, Deep BSDE \cite{han2018solving,mastrolia2025optimal}, and PINNs \cite{raissi2019physics,gennaro2024delegated,guo2026actor}.
We give the main lines of the algorithm we used in Appendix \ref{app:alg_conv} together with illustrating the convergence and reduction of the loss functions associated with our simulations in Figure \ref{fig:Conv_Graph}.

We now turn to the numerical illustration of our results. 
Figure \ref{fig:One_Path_Graphs} shows one simulated path of the state process $(S,I,P)$, the jumps processes and optimizers for both agent and principal while Figure \ref{fig:32_Path_Graphs} reports 32 simulated paths of these processes summarized by their mean and a 90\% confidence band. The numerical results in both Figure \ref{fig:One_Path_Graphs} (a) and Figure \ref{fig:32_Path_Graphs} (a) show that the price process under the contract is systematically higher than in the no-contract benchmark. The mechanism is visible in the SIR dynamics: under the contract, the susceptible population declines more rapidly, while the infected population rises more slowly and remains lower. Figure \ref{fig:One_Path_Graphs} (b) and Figure \ref{fig:32_Path_Graphs} (b) illustrate the impact of the contract on both external and internal jumps. We see that there is no impact on the external jumps since it is only induced by the attacker and the defender (agent) is not controlling it, the difference observed in the external jumps only comes from the fact that the price drops in percentage (see \eqref{eq:P-Lhop-random}). Regarding the internal losses induced by the attack, we notice a significant reduction induced by the optimal action of the agent. In our reduced SIR representation, this means that the contract induces the subsidiary to move devices more quickly from the susceptible class into the recovered class, thereby limiting the propagation of infection, preserving project value while preventing infection outbreak and future losses.

This is also confirmed by the incentive simulations shown in Figures \ref{fig:One_Path_Graphs} and \ref{fig:32_Path_Graphs} (c) and (d)  below. The contract induces substantially stronger effort, especially early in the horizon, when preventive action is most valuable. This front-loading of protection is economically natural: when a large susceptible population is still exposed, marginal cybersecurity effort has the largest effect on the future infection path and on the associated jump exposure. The dominant component is the continuous loading $Z_2$, which is the part of the contract that directly incentives preventive effort in the susceptible-state direction. By contrast, the external jump-loading is fixed at zero in the final calibration, reflecting the fact that external attacks are not directly controlled by the subsidiary. The internal jump-loading $U^i$ is active and negative, which means that the agent's continuation utility falls after an internal cyber incident. Since the subsidiary's effort lowers the internal jump intensity, this jump-contingent penalty provides an additional incentive to protect the system. Numerically, the optimal contract is therefore prevention-dominant: it works mainly through a strong continuous incentive component together with an internal jump penalty, rather than through insurance against realized losses.\\

Figure \ref{fig:principal_values} reports the principal's value with and without the contract, as well as the contract gain for different initial conditions $(s_0=s,p_0=s)$ and $i=1-s_0$ assuming no recovered device at the beginning. The gain is positive on the plotted state space and is largest in regions with a relatively large susceptible population. This pattern is consistent with the preventive nature of the contract: when more devices remain vulnerable, stronger effort can materially alter the future infection dynamics, so the contract has more leverage. Regarding the absolute contract gain (b) we observe that the gain is more significant if $p_0$ and $s_0$ are higher. It is explained by the fact that if $s_0$ is large enough then the agent has more margin to protect devices and so make the system more resilient to cyber attacks. As the price higher at the beginning induced a higher value due to the price model chosen and the reduction of L-hop risk.\\

Finally, we complement the heatmaps with a compact fixed-state sensitivity analysis by fixing $(p_0,s_0,i_0)$ and perturbing one category of cyber-risk parameters at a time, see Figure \ref{fig:Sen_Ana}. The main object of interest is the contract gain \textit{i.e.} the increase in the principal's value from access to the contract relative to the no-contract benchmark. Across all tested perturbations, the contract gain remains positive, showing that the contracting mechanism is robustly valuable at the selected state. The sensitivity experiments also clarify how the value of contracting depends on the nature of cyber risk. When jump sizes are increased, the gain declines moderately, which indicates that larger cyber incidents create more realized loss that cannot be fully eliminated once the jump occurs. When attack intensities are increased, the gain remains positive but can again decline in the most severe regimes, reflecting the fact that very frequent attacks generate common losses for both regimes. By contrast, increasing the infection sensitivity of internal jump arrivals preserves a relatively high contract gain. This is the scenario most closely aligned with the mechanism of the model: the contract raises preventive effort, lowers the infection path, and thereby reduces the endogenous component of internal cyber-jump risk. Overall, the sensitivity analysis suggests that the contract is most valuable in moderate-to-high cyber-risk environments where prevention can still materially alter the infection dynamics. In extremely severe regimes, the contract remains robust and ensures a higher system resiliency.

}

\newpage

\bibliographystyle{apalike}
\bibliography{siam_review/siam_Mastrolia_Yan/references}

\appendix

\section{Spaces}\label{app:spaces}
Let $ \mathbb X := (\mathcal X_s)_{t \leq s \leq T} $ denote an arbitrary filtration on $ (\Omega, \mathcal{F}_T) $, and let $ \mathbb{P} $ be an arbitrary element in $ \mathcal{P}(t, \omega) $. We follow the notations of spaces in \cite{vuca,possamai2018stochastic,denis:hal-04822047}. 
\begin{itemize}
\item \textbf{The spaces $\mathbb{L}^{p,\kappa}_{t,x}$.} For each $ p \geq \kappa \geq 1 $, we define $\mathbb L^p_{t,\omega}(\mathbb X)$  (resp.  $\mathbb L^p_{t,\omega}(\mathbb X, \mathbb{P}))$ denotes the space of all $\mathcal{X}_T-$measurable random variables $\xi$ such that
$$
\|\xi\|_{\mathbb L^p_{t,\omega}} := \sup_{\mathbb{P} \in \mathcal{P}(t,\omega)} \left( \mathbb{E}^{\mathbb{P}}[|\xi|^p] \right)^{1/p} < +\infty,$$
and respectively
$$\|\xi\|_{\mathbb L^p_{t,\omega}(\mathbb{P})} := \left( \mathbb{E}^{\mathbb{P}}[|\xi|^p] \right)^{1/p} < +\infty.
$$
We set
$$
\mathbb L^{p,\kappa}_{t,\omega}(\mathbb X) := \left\{ \xi \in \mathbb L^p_{t,\omega}(\mathbb X) : \|\xi\|_{\mathbb L^{p,\kappa}_{t,\omega}} < \infty \right\},
$$
where the norm is given by
$$
\|\xi\|_{\mathbb L^{p,\kappa}_{t,\omega}} := \sup_{\mathbb{P} \in \mathcal{P}(t,\omega)} 
\left( 
\mathbb{E}^{\mathbb{P}} \left[ 
\operatorname*{esssup}_{t \leq s \leq T} 
\left( \mathbb{E}^{\mathbb{P}}_{t,\omega,\mathcal{X}_s^+} \left[ |\xi|^\kappa \right] \right)^{\frac{p}{\kappa}} 
\right] 
\right)^{\frac{1}{p}}.
$$
\item 
\textbf{The spaces $\mathbb{H}^{p}_{t,x}(\mathbb X,\mathbb P)$.} We say $Z$ is in $\mathbb{H}^{p}_{t,x}(X,\mathbb P)$ if $Z$ is an $X$-predictable, $\mathbb{R}^d-$valued process satisfying
$$
\|Z\|_{\mathbb{H}^{p}_{t,x}(\mathbb X,\mathbb P)}^p := \mathbb{E}^{\mathbb P}\Biggl[\biggl(\int_{t}^T \|\sigma_s^{\tfrac12}Z_s\|^2ds\biggr)^{\tfrac{p}{2}}\Biggr] < +\infty.
$$
We then define
$$
\mathbb{H}^{p}_{t,x}(X,\mathcal{P}) := \Bigl\{Z:\sup_{\mathbb P \in \mathcal{P}(t,x)} \|Z\|_{\mathbb{H}^{p}_{t,x}(\mathbb X,\mathbb P)} < +\infty\Bigr\}.
$$

\item
\textbf{The spaces $\mathbb{S}^{p}_{t,x}(\mathbb X,\mathbb P)$.} We say $Y$ is in $\mathbb{S}^{p}_{t,x}(\mathbb X,\mathbb P)$ if $Y$ is an $\mathbb X$-progressively measurable, real-valued process satisfying
$$
\|Y\|_{\mathbb{S}^{p}_{t,x}(\mathbb X,\mathbb P)}^p := \mathbb{E}^P\Bigl[\sup_{s\in [t,T]} |Y_s|^p\Bigr] < +\infty.
$$
We then define
$$
\mathbb{S}^{p}_{t,x}(\mathbb X,\mathcal{P}) := \Bigl\{Y:\sup_{\mathbb P \in \mathcal{P}(t,x)} \|Y\|_{\mathbb{S}^{p}_{t,x}(\mathbb X,\mathbb P)} < +\infty\Bigr\}.
$$

\item \textbf{ The space $\mathbb J^p_{t,\omega}(\mathbb X)$.} $\mathbb J^p_{t,\omega}(\mathbb X)$ (resp. $ \mathbb J^p_{t,\omega}(\mathbb X, \mathbb{P}) $) denotes the space of all $\mathbb X$-predictable functions $U$ such that
$$
\|U\|_{\mathbb J^p_{t,\omega}(\mathbb X)}\! := \!\sup_{\mathbb{P} \in \mathcal{P}(t,\omega)} 
\!\left( 
\mathbb{E}^{\mathbb{P}} \left[ 
\left( \int_t^T \int_{\mathbb R\setminus\{0\}} \|U_s(\chi)\|^2  \lambda_s^0(d\chi)  ds \right)^{p/2} 
\right] 
\right)^{1/p} < +\infty,
$$
resp.
$$
\|U\|_{\mathbb J^p_{t,\omega}(\mathbb X,\mathbb{P})} := 
\left( 
\mathbb{E}^{\mathbb{P}} \left[ 
\left( \int_t^T \int_{\mathbb R\setminus\{0\}} \|U_s(\chi)\|^2  \lambda_s^0(d\chi)  ds \right)^{p/2} 
\right] 
\right)^{1/p} < +\infty.
$$
\item
\textbf{The spaces $\mathbb{K}^{p}_{t,x}(\mathbb X,\mathbb P)$.} We say $K$ is in $\mathbb{K}^{p}_{t,x}(\mathbb X,\mathbb P)$ if $K$ is an $\mathbb X$-optional, real-valued process with $\mathbb P-$a.s. càdlàg, non-decreasing paths on $[t,T]$, $K_t=0$ $\mathbb P-$a.s., and
$$
\|K\|_{\mathbb{K}^{p}_{t,x}(\mathbb X,\mathbb P)}^p := \mathbb{E}^P\bigl[|K_T|^p\bigr] < +\infty.
$$
We denote by $\mathbb{K}^{p}_{t,x}\bigl(\mathbb X,\mathcal{P}\bigr)$ the set of all families $(K^P)_{P\in \mathcal{P}(t,x)}$ such that $K^P \in \mathbb{K}^{p}_{t,x}(\mathbb X,\mathbb P)$ for every $\mathbb P \in \mathcal{P}(t,x)$ and
$$
\sup_{\mathbb P\in \mathcal{P}(t,x)} \|K^{\mathbb P}\|_{\mathbb{K}^{p}_{t,x}(\mathbb X,\mathbb P)} < +\infty.
$$

\item
\textbf{The spaces $\mathcal L^{p}_\nu$.}

We define $\mathcal L^{p}_\nu$ as the set of Borel measurable functions $u: \mathbb{R}^* \to \mathbb{R}^m$ satisfying
$$
\|u\|_{p,\nu} := \int_{\mathbb{R}^*} \|u(\chi)\|^p\nu(d\chi) < +\infty.
$$

\item
{\textbf{The space \(C^{1,2}_{\mathrm{pol}}\).}
For the Perron argument, the state variable is \((x,y)\in\mathbb R^n\times\mathbb R\).
We write
$$
D_{x,y}\varphi
:=
\big(D_x\varphi,\partial_y\varphi\big),
\qquad
D^2_{x,y}\varphi
:=
\begin{pmatrix}
D^2_{xx}\varphi & D_{xy}\varphi\\
D_{yx}\varphi & \partial^2_{yy}\varphi
\end{pmatrix}.
$$
For \(m\ge 1\), we define \(C^{1,2}_{\mathrm{pol},m}\big([0,T]\times
\mathbb R^n\times\mathbb R\big)\) as the set of all functions
$$
\varphi\in C\big([0,T]\times\mathbb R^n\times\mathbb R\big)
\cap
C^{1,2}\big([0,T)\times\mathbb R^n\times\mathbb R\big)
$$
such that
$$
\|\varphi\|_{1,2;m}
:=
\sup_{(t,x,y)\in[0,T]\times\mathbb R^n\times\mathbb R}
\frac{|\varphi(t,x,y)|}{1+|x|^m+|y|^m}
$$
$$
\quad
+
\sup_{(t,x,y)\in[0,T)\times\mathbb R^n\times\mathbb R}
\frac{
|\partial_t\varphi(t,x,y)|
+
|D_{x,y}\varphi(t,x,y)|
+
|D^2_{x,y}\varphi(t,x,y)|
}{
1+|x|^m+|y|^m
}
<\infty.
$$
We then set
$$
C^{1,2}_{\mathrm{pol}}
\big([0,T]\times\mathbb R^n\times\mathbb R\big)
:=
\bigcup_{m\ge 1}
C^{1,2}_{\mathrm{pol},m}
\big([0,T]\times\mathbb R^n\times\mathbb R\big).
$$

Whenever there is no ambiguity, we simply write \(C^{1,2}_{\mathrm{pol}}\).
Thus a function in \(C^{1,2}_{\mathrm{pol}}\) is continuous up to the terminal time,
is \(C^1\) in time and \(C^2\) in the space variables on \([0,T)\), and its value,
time derivative, first spatial derivatives, and second spatial derivatives have
polynomial growth.
}

\end{itemize}

\section{Proof of Theorem~\ref{thm:agent-2BSDEJ}}

We follow the scheme in \cite{vuca}.
We first prove that \eqref{eq:ver} holds with a characterization of the optimal effort of the Agent as a maximizer of the 2BSDEJ \eqref{eq:2BSDEJ-agent}. The proof is divided into five steps.

\textbf{Step 1: BSDEJ and 2BSDJ.} 
For every $\bigl(\alpha,\eta\bigr) \in \mathcal{A} \times \mathcal{H}(\hat{\sigma})$, denote by 
$\bigl(Y^{\alpha,\eta}, Z^{\alpha,\eta}, U^{\alpha,\eta},K^{\alpha,\eta}\bigr)$
the solution of the following controlled 2BSDEJ in the sense of Definition \ref{2bsde:def} and where the wellposedness is deduced from \cite{denis:hal-04822047}.
\begin{equation}\label{eq:2BSDEJ-fixed-a-eta}
\begin{aligned}
Y^{\alpha,\eta}_t
&=U^A(\xi)+F^A(X_T)
+\int_t^T G\big(s,X_s,Y^{\alpha,\eta}_s,Z^{\alpha,\eta}_s,U^{\alpha,\eta}_s;\alpha_s,\eta_s\big)ds\\
& -\int_t^T Z^{a}_s\cdot dX^{c,\mathbb P}_s
      -\int_t^T\int_{\mathbb R^n\setminus\{0\}} U^{\alpha,\eta}_s(\chi)(\mu_X-\Lambda^0)(ds,d\chi)
      -\int_t^T dK^{\alpha,\eta}_s,
\mathcal P-q.s.
\end{aligned}
\end{equation}

Note in particular, see \cite[Section 2.5]{denis:hal-04822047} and \cite[Theorem 4.2]{possamai2018stochastic} that 
\begin{equation}\label{eq:rep-fixed-a-h}
    Y_0^{\alpha,\eta}
  =
  {\essinf_{\mathbb{P}'\in\mathcal{P}[\mathbb{P},\mathbb{F}^+,0]}}^\mathbb{P}
  \mathcal Y_0^{\mathbb{P}',\alpha,\eta},
  \quad
  \mathbb{P}\text{-a.s.\ for every } \mathbb{P}\in\mathcal{P}.
\end{equation}
where for any $\mathbb P\in\mathcal P$ the tuple 
$\bigl(\mathcal Y^{\mathbb{P},u,\alpha}_t,\mathcal Z^{\mathbb{P},u,\alpha}_t,\mathcal U_t^{\mathbb{P},u,\alpha}\bigr)$
is the solution of the following (well-posed) linear BSDEJ, see for example \cite{papapantoleon2018existence}
\begin{equation}\label{eq:BSDEJ-fixed-a-eta}
\begin{aligned}
\mathcal Y_t^{\mathbb P;a,\eta}
&=U^A(\xi)+F^A(X_T)\\
&\quad+\int_t^T G\big(s,X_s,\mathcal Y_s^{\mathbb P;a,\eta},\mathcal Z_s^{\mathbb P;a,\eta},\mathcal U_s^{\mathbb P;a,\eta};a_s,\eta_s\big)ds\\
&\quad-\int_t^T \mathcal Z_s^{\mathbb P;a,\eta}\cdot dX^{c,\mathbb P}_s
      -\int_t^T\int_{\mathbb R^n\setminus\{0\}} \mathcal U_s^{\mathbb P;a,\eta}(\chi)(\mu_X-\Lambda^0)(ds,d\chi),
\qquad \mathbb P\text{-a.s.}
\end{aligned}
\end{equation}

Similarly, consider also for each $a\in\mathfrak A$, for all $\mathbb P\in\mathcal P$,
\begin{equation}\label{eq:2BSDEJ-fixed-a-muX}
\begin{aligned}
Y^{a}_t
&=U^A(\xi)+F^A(X_T)
+\int_t^T \inf_{\eta\in \mathcal H(s,X_s,\widehat\sigma_s)}G\big(s,X_s,Y^{a}_s,Z^{a}_s,U^{a}_s;\ \alpha,\eta)ds\\
& -\int_t^T Z^{\alpha}_s\cdot dX^{c,\mathbb P}_s
      -\int_t^T\int_{\mathbb R^n\setminus\{0\}} U^{\alpha}_s(\chi)(\mu_X-\Lambda^0)(ds,d\chi)
      -\int_t^T dK^{\alpha}_s,
\quad \mathbb P\text{-a.s., }
\end{aligned}
\end{equation}
By the standard 2BSDEJ representation (upper envelope of single-prior BSDEJs on the set of continuations),
\begin{equation}\label{eq:rep-fixed-a}
Y^{a}_0
=
\essinf_{\mathbb P'\in\mathcal P[\mathbb P,\mathcal F^+_0,0]}
\mathcal Y^{\mathbb P';a,\eta^{\mathbb P'}}_0,
\quad\text{where}\quad
\eta^{\mathbb P'}_s\in \arg\min_{\eta\in \mathcal H(s,X_s,\widehat\sigma_s)} G(\cdot;a_s,\eta).
\end{equation}

\smallskip
\noindent
\textbf{Step 2 (comparison across $a$ and reconstruction of $G^*$).}
From comparison theorem for the BSDEJ, we deduce that  $\mathcal Y_0^{\mathbb P,\alpha}\leq \mathcal Y_0^{\mathbb P,\alpha,\eta},$ for any $\mathbb P\in \mathcal P$ and the equality hold for $\eta$ optimizing the infimum. 
Therefore, from the representation \eqref{eq:rep-fixed-a-h} and \eqref{eq:rep-fixed-a} we deduce that
\begin{equation}\label{eq:rep}
  Y_0 
  = \esssup_{\alpha \in \mathcal{A}} Y_0^\alpha
  = \esssup_{\alpha \in \mathcal{A}} 
    \essinf_{\eta \in \mathcal{H}(\hat{\sigma}^2)}Y_0^{\alpha,\eta},
  \quad \mathbb{P}\text{-a.s.\ for every } \mathbb{P}\in \mathcal{P}.
\end{equation}

\textbf{Step 3: linearization and value function.} The generator $G$ is linear in $y,z,u$. By using standard linearization tools for BSDEJ, see for example \cite{quenez2013bsdes} we get $$
  \mathcal Y_0^{\mathbb{P},\alpha,\eta }
  =
  \mathbb E^{\mathbb P}\left[\mathcal K_{0,T}\Big(U^A(\xi) + F^A(\mathbf{X}_T)\Big) - \int_0^T \mathcal K_{0,s}C^A(s,\mathbf{X}_s,\alpha_s) ds\right],
  \quad
  \mathbb{P}\text{-a.s.}, \mathbb P\in \mathcal P_0.
$$

\textbf{Step 4: characterization of the value function.}
From the previous steps, it follows that $\mathbb{P}^{\alpha,\eta} \in \mathcal{P}^{\alpha}$ and 
$\mathbb{P}$-a.s.\ for every $\mathbb{P}\in\mathcal{P}$:
\begin{equation*}
    \begin{split}
        Y_0
        &=
        {\esssup_{\alpha\in\mathcal{A}}}
        {\essinf_{\eta\in\mathcal{H}(\hat\sigma^2)}}
        {\essinf_{\mathbb{P}'\in\mathcal{P}[\mathbb{P},\mathbb{F}^+,0]}}
       \mathbb E^{\mathbb P}\left[\mathcal K_{0,T}\Big(\!U^A(\xi)\!+ \!F^A({X}_T)\!\Big)\! - \!\int_0^T \mathcal K_{0,s}C^A(s,\!{X}_s,\!\alpha_s) ds\right]\\
        & ={\esssup_{\alpha\in\mathcal{A}}}{\essinf_{(\mathbb{P}',\eta)\in\mathcal{H}^\alpha[\mathbb{P},\mathbb{F}^+,0]}}\mathbb E^{\mathbb P}\left[\mathcal K_{0,T}\Big(U^A(\xi)\! + \!F^A({X}_T)\Big) \!-\! \int_0^T \mathcal K_{0,s}C^A(s,{X}_s,\alpha_s) ds\right].
    \end{split}
\end{equation*}
  
The characterization \eqref{eq:ver} then follows by similar arguments to those used in the proofs of Lemma 3.5 and Theorem 5.2 of \cite{possamai2018stochastic}.

\textbf{Step 5: optimizers.}
We now turn to the second part of the theorem, where the characterization of an optimal triplet $(\alpha,\eta,\mathbb{P})$ for the optimization problem \eqref{eq:ver} is shown.  From the previous steps, it is clear that a control $\bigl(\hat\alpha,\eta^\star,\mathbb{P}^\star\bigr)$ is optimal if and only if it attains all the essential suprema and infima above.  In particular, the infimum in \eqref{eq:rep-fixed-a-h} is attained under conditions $({ii})$, and equality \eqref{eq:rep} holds if $\bigl(\hat\alpha,\eta^*\bigr)$ satisfy ${(i)}$.

{ \section{Proof of Theorem~\ref{thm:viscosity}}
\label{app:viscosity_proof}

\begin{lemma}[Elementary approximation with jumps]
\label{lem:vuca-K-regularization-jumps}
Fix bounded local controls $(z,u,\gamma)$.
Set
\begin{equation*}\label{eq:C-Kdot-vuca}
\begin{split}
\dot K_s^{z,u,\gamma}
&:=
G^\star(s,X_{s-},Y_{s-},z,u;\widehat\sigma_s)
+\frac12\Tr(\widehat\sigma_s\gamma)
-H(s,X_{s-},Y_{s-},z,u;\gamma),\\
K_s^{z,u,\gamma}
&:=
\int_t^s\dot K_r^{z,u,\gamma}dr .
\end{split}
\end{equation*}
Then $K^{z,u,\gamma}$ is nondecreasing. Moreover, for every bounded predictable
process $\psi$ and every $\varepsilon>0$, there exists a nonnegative elementary
predictable density $k^p$, piecewise constant,
such that, for all stopping times $\theta\le\rho$ where $\rho$ is a localization time ensuring that $X,Y$ are bounded
\begin{equation}\label{eq:C-K-approx}
\left|
\int_\theta^\rho
\psi_s\big(\dot K_s^{z,u,\gamma}-k^p_s\big)ds
\right|
\le
\varepsilon(\rho-\theta),
\qquad
\mathfrak P^{(\widehat\alpha)}\text{-q.s.}
\end{equation}
\end{lemma}

\begin{proof}
The nonnegativity of $\dot K^{z,u,\gamma}$ follows directly from the definition
of the Hamiltonian: for $\Sigma=\widehat\sigma_s$,
$$
H(s,X_{s-},Y_{s-},z,u;\gamma)
\le
G^\star(s,X_{s-},Y_{s-},z,u;\widehat\sigma_s)
+\frac12\Tr(\widehat\sigma_s\gamma).
$$
The integrand $\dot K^{z,u,\gamma}$ is bounded and
predictable on $[0,\rho]$. Hence it can be approximated, path by path in the Lebesgue integral,
by nonnegative elementary predictable processes on finite stopping-rule grids.
This gives~\eqref{eq:C-K-approx}. This is exactly the approximation argument of
Lemma D.1 in \cite{vuca}. The presence of jumps does not change the argument,
because the approximation concerns the absolutely continuous Lebesgue integral
in time; jumps only enter through the predictable càdlàg arguments
$(X_{s-},Y_{s-})$ and the fixed jump loading $u$.
\end{proof}

\begin{lemma}[Coercive localization for jump]
\label{lem:coercive-localization}
Let $\varphi\in C^{1,2}_{\rm pol}$ such that $v^{-}-\varphi$ reaches a strict local minimum at
$(t_0,x_0,y_0)\in[0,T)\times\mathbb R^n\times\mathbb R$, and suppose that
$$
\partial_t\varphi(t_0,x_0,y_0)
+
\mathcal Q^\star[\varphi](t_0,x_0,y_0)>0.
$$
Then there exist $q\ge1$, $\iota>0$, bounded neighborhoods
$\mathcal O_0\subset\mathcal O$ of $(t_0,x_0,y_0)$, as for example a ball in $\mathbb R^{n+2}$ centered in $(t_0,x_0,y_0)$ with radius some positive integer $m$ chosen with
$\overline{\mathcal O}\subset[0,T)\times\mathbb R^n\times\mathbb R$, a constant
$\kappa>0$, a stochastic subsolution $w\in\mathcal V^{-}$, a constant $a>0$,
and localized controls $(\widehat z,\widehat u,\widehat\gamma)$ such that, with
$$
\ell(x,y):=|x-x_0|^{2q}+|y-y_0|^{2q},
\qquad
\widetilde\varphi(t,x,y):=\varphi(t,x,y)-\iota\ell(x,y),
$$
we have
\begin{equation}\label{eq:C-vminus-strict-localized}
\partial_t\widetilde\varphi
+
\mathcal Q^{\widehat z,\widehat u,\widehat\gamma,\eta}
[\widetilde\varphi]
\ge 3a
\quad\text{on }\mathcal O,
\qquad \forall\eta\in H,
\end{equation}
and
\begin{equation}\label{eq:C-vminus-global-order}
\widetilde\varphi+\kappa\le w
\quad\text{on }\mathcal O\setminus\mathcal O_0.
\end{equation} 

Symmetrically, if $\varphi\in C^{1,2}_{\rm pol}$ such that $v^{+}-\varphi$ reaches a strict maximum at the point
$(t_0,x_0,y_0)$ and
$$
\partial_t\varphi(t_0,x_0,y_0)
+
\mathcal Q^\star[\varphi](t_0,x_0,y_0)<0,
$$
then there exist $q\ge1$, $\iota>0$, bounded neighborhoods
$\mathcal O_0\subset\mathcal O$ of $(t_0,x_0,y_0)$, chosen with
$\overline{\mathcal O}\subset[0,T)\times\mathbb R^n\times\mathbb R$, $\kappa>0$,
$w\in\mathcal V^{+}$, $a>0$, and a measurable selector
$$
\widehat\eta=\widehat\eta(t,x,y,z,u,\gamma)\in H
$$
such that, with
$$
\ell(x,y):=|x-x_0|^{2q}+|y-y_0|^{2q},
\qquad
\widetilde\varphi(t,x,y):=\varphi(t,x,y)+\iota\ell(x,y),
$$
we have, for every localized $(z,u,\gamma)$,
\begin{equation}\label{eq:C-vplus-strict-localized}
\partial_t\widetilde\varphi
+
\mathcal Q^{z,u,\gamma,\widehat\eta(t,x,y,z,u,\gamma)}
[\widetilde\varphi]
\le -3a
\quad\text{on }\mathcal O,
\end{equation}
and
\begin{equation}\label{eq:C-vplus-global-order}
w\le\widetilde\varphi-\kappa
\quad\text{on }\mathcal O\setminus\mathcal O_0,
\end{equation}
with the same inequality holding at every possible post-jump exit point
$$
\big(
s,
x+\beta(s,x,e),
y+u(\beta(s,x,e))
\big)\notin\mathcal O_0
$$
arising from $(s,x,y)\in\mathcal O_0$, $e\in E$, and localized $u$.
\end{lemma}
\begin{remark}\label{rem:jumpvuca}
    Note that the $\iota \ell$ terms in the lemma above is specific to the jump case compared with \cite{vuca}. It ensures that the inequality \eqref{eq:C-vminus-global-order} remains true regarding any jumps, see \cite[Appendix B]{BayraktarLi2016}.
\end{remark}
\begin{proof}[Proof of Lemma \ref{lem:coercive-localization}]
By choosing $q$ larger than the polynomial-growth exponent of the payoff and
semisolutions, the coercive term $\iota \ell$ controls the possible post-jump values. The stability of stochastic subsolutions
under finite maxima follows from the elementary conditioning and pasting argument
used in Lemmas 3.7--3.8 of \cite{perronsirbu}; the same jump-target pasting idea is
used in Lemma 3.1 of \cite{BayraktarLi2016}. Therefore the finite maximum belongs
to $\mathcal V^{-}$ and yields the desired $w$. The $v^{+}$ case is identical, with finite minima replacing finite maxima.
\end{proof}

\begin{proof}[Proof of Theorem~\ref{thm:viscosity}]
The proof follows the stochastic Perron argument of \cite{vuca}. We only record
the changes caused by jumps: the elementary Principal control is now
$(Z,U,K)$, the It\^o formula is the It\^o--Lévy formula together with Lemma \ref{lem:coercive-localization}. We first prove that $v^{-}$ is a lower semicontinuous viscosity supersolution and that
$v^{+}$ is an upper semicontinuous viscosity subsolution of~\eqref{eq:HJBI-ref}
on $[0,T)\times\mathbb R^n\times\mathbb R$. We also prove
$$
(v^{-})_\ast(T,x,y)\ge g(x,y),
\qquad
(v^{+})^\ast(T,x,y)\le g(x,y).
$$

\medskip
\noindent
\textbf{Step 1. It\^o formula with jump and localization.}
Let $\phi\in C^{1,2}_{\rm pol}$. Fix local controls $(z,u,\gamma)$, let $k^p$
be the elementary approximation of $\dot K^{z,u,\gamma}$ given by
Lemma~\ref{lem:vuca-K-regularization-jumps}, and fix an admissible continuation
pair $(\mathbb P,\eta)$. For stopping times $\theta\le\rho$ is a localization time ensuring that $X,Y$ are bounded (as for example a ball centered in $(x_0,y_0)$ with radius $m>0$), the It\^o-Lévy formula gives
\begin{equation}\label{eq:C-jump-ito-cost}
\begin{aligned}
\phi(\rho,X_\rho,Y_\rho)
&=
\phi(\theta,X_\theta,Y_\theta)+
\int_\theta^\rho
C^P(s,X_{s-},\widehat\alpha_s,\eta_s)\,ds +M_\rho-M_\theta \\
&
+
\int_\theta^\rho
\left(
\partial_t\phi
+
\mathcal Q^{z,u,\gamma,\eta}[\phi]
\right)(s,X_{s-},Y_{s-})\,ds\\
&+
\int_\theta^\rho
\partial_y\phi(s,X_{s-},Y_{s-})
\big(k_s^p-\dot K_s^{z,u,\gamma}\big)\,ds.
\end{aligned}
\end{equation}
Here $M$ is a local martingale. Its continuous part is
$$
\int_\theta^\rho
\Big(
D_x\phi(s,X_{s-},Y_{s-})
+
\partial_y\phi(s,X_{s-},Y_{s-})z
\Big)^\top
\sigma(s,X_{s-},\eta_s)
\,dW_s^{(\widehat\alpha,\eta)} ,
$$
and its compensated jump part is
$$
\int_\theta^\rho\int_E
\Big[
\phi\big(
s,
X_{s-}+\beta(s,X_{s-},e),
Y_{s-}+u(\beta(s,X_{s-},e))
\big)
-
\phi(s,X_{s-},Y_{s-})
\Big]
\widetilde\mu^{(\widehat\alpha,\eta)}(ds,de).
$$
The compensator of this jump martingale is exactly the nonlocal term in
$\mathcal L^{z,u,\gamma,\eta}\phi$. By Assumption~\ref{ass:perron-standing}(ii),
after the standard localization the martingale is a true martingale.

\medskip
\noindent
\textbf{Step 2. $v^{-}$ is a viscosity supersolution on $[0,T)$.}

\smallskip
\noindent
\emph{a. Setting the contradiction.}
Let $\varphi\in C^{1,2}_{\rm pol}$ be such that $v^{-}-\varphi$ has a strict
local minimum equal to zero at $(t_0,x_0,y_0)\in[0,T)\times\mathbb R^n\times
\mathbb R$. We assume by contradiction that
\begin{equation*}\label{eq:C-vminus-contradiction}
\partial_t\varphi(t_0,x_0,y_0)
+
\mathcal Q^\star[\varphi](t_0,x_0,y_0)>0 .
\end{equation*}
Using Lemma~\ref{lem:coercive-localization}, and shrinking the time component
of the patch if necessary so that $\mathcal O_0\subset[0,T)\times\mathbb R^n
\times\mathbb R$, there exist $\widetilde\varphi$, $\mathcal O_0\subset
\mathcal O$, $\kappa>0$, $a>0$, localized controls
$(\widehat z,\widehat u,\widehat\gamma)$, and $w\in\mathcal V^{-}$ such that
\eqref{eq:C-vminus-strict-localized}--\eqref{eq:C-vminus-global-order} hold.

Define
\begin{equation*}\label{eq:C-vminus-patch}
w^\kappa(t,x,y)
:=
\begin{cases}
\max\{w(t,x,y),\widetilde\varphi(t,x,y)+\kappa\},
& (t,x,y)\in\mathcal O_0,\\ [1mm]
w(t,x,y),
& (t,x,y)\notin\mathcal O_0 .
\end{cases}
\end{equation*}
By the global ordering \eqref{eq:C-vminus-global-order}, $w^\kappa$ is
continuous. Moreover,
$$
w^\kappa(t_0,x_0,y_0)
\ge
\widetilde\varphi(t_0,x_0,y_0)+\kappa
>
v^{-}(t_0,x_0,y_0).
$$
Thus it remains to prove that $w^\kappa\in\mathcal V^{-}$.

\smallskip
\noindent
\emph{b. Building the elementary strategy and verifying property $(ii-)$.}
Fix $t\in[0,T]$ and a stopping rule $\tau$. Since $w\in\mathcal V^{-}$, the
definition of stochastic subsolution gives elementary continuation controls for
$w$ at $\tau$ and, after the exit time below, at $\theta$.

Let $\tau'$ be the stopping time induced by $\tau$ along the controlled path and
set
$$
A^{-}
:=
\left\{
\Xi_{\tau'}\in\mathcal O_0,
\quad
\widetilde\varphi(\Xi_{\tau'})+\kappa>w(\Xi_{\tau'})
\right\}
\in\mathcal F_{\tau'}^t .
$$
Let
$$
\theta
:=
\inf\{s\ge\tau:\Xi_s\notin\mathcal O_0\}\wedge T .
$$
On $A^{-}$ we use the local controls
$$
(\widehat z,\widehat u,\widehat k^p)
\quad\text{on }(\tau,\theta],
$$
where $\widehat k^p$ is given by
Lemma~\ref{lem:vuca-K-regularization-jumps} applied to
$\dot K^{\widehat z,\widehat u,\widehat\gamma}$ with
$\psi=\partial_y\widetilde\varphi$ and error smaller than the strict margin
$a$. After $\theta$ we use the continuation control associated with $w$. On
$(A^{-})^c$ we use the continuation control associated with $w$ directly.
By Assumption~\ref{ass:perron-standing}(iii), the resulting pasted control
belongs to $\mathfrak K_{Y_0}(t,\tau)$.

Now fix an arbitrary prior elementary Principal control, an arbitrary
continuation Nature pair $(\mathbb P,\eta)$, and an arbitrary stopping rule
$\rho$ with $\tau\le\rho\le T$. Denote by $\rho'$ and $\theta'$ the induced
stopping times. On $A^{-}$, applying \eqref{eq:C-jump-ito-cost} to
$\widetilde\varphi+\kappa$ on $[\tau',\rho'\wedge\theta']$, and using
\eqref{eq:C-vminus-strict-localized} together with
\eqref{eq:C-K-approx}, shows that the cost-adjusted process
$$
\widetilde\varphi(\Xi_s)+\kappa
-
\int_{\tau'}^s
C^P(r,X_{r-},\widehat\alpha_r,\eta_r)\,dr
$$
is a submartingale on the stopped local branch. Hence, on $A^{-}$,
\begin{equation*}\label{eq:C-vminus-local-submartingale}
\widetilde\varphi(\Xi_{\tau'})+\kappa
\le
\mathbb E^{\mathbb P}
\left[
\widetilde\varphi(\Xi_{\rho'\wedge\theta'})+\kappa
-
\int_{\tau'}^{\rho'\wedge\theta'}
C^P(s,X_{s-},\widehat\alpha_s,\eta_s)\,ds
\middle|
\mathcal F_{\tau'}^t
\right].
\end{equation*}

Inside $\mathcal O_0$ we have
$w^\kappa\ge\widetilde\varphi+\kappa$. If the process exits
$\mathcal O_0$ by a jump, the post-jump point belongs to $\mathcal O_0^c$ and
the global ordering \eqref{eq:C-vminus-global-order} gives
$$
w^\kappa=w\ge\widetilde\varphi+\kappa
\quad\text{on }\mathcal O_0^c .
$$
Thus the preceding estimate is stable under jump overshoots and yields
$$
w^\kappa(\Xi_{\tau'})
\le
\mathbb E^{\mathbb P}
\left[
w^\kappa(\Xi_{\rho'\wedge\theta'})
-
\int_{\tau'}^{\rho'\wedge\theta'}
C^P(s,X_{s-},\widehat\alpha_s,\eta_s)\,ds
\middle|
\mathcal F_{\tau'}^t
\right]
\quad\text{on }A^{-}.
$$
On $(A^{-})^c$, the same inequality follows from the stochastic subsolution
property of $w$, using $w^\kappa\ge w$. After $\theta'$, the control follows the
continuation control for $w$; since $w^\kappa=w$ at the switching time and
$w^\kappa\ge w$ everywhere, the tower property gives
$$
w^\kappa(\Xi_{\tau'})
\le
\mathbb E^{\mathbb P}
\left[
w^\kappa(\Xi_{\rho'})
-
\int_{\tau'}^{\rho'}
C^P(s,X_{s-},\widehat\alpha_s,\eta_s)\,ds
\middle|
\mathcal F_{\tau'}^t
\right].
$$
Therefore $w^\kappa\in\mathcal V^{-}$, contradicting the definition of $v^{-}$.
Consequently,
$$
\partial_t\varphi(t_0,x_0,y_0)
+
\mathcal Q^\star[\varphi](t_0,x_0,y_0)
\le0.
$$
Hence $v^{-}$ is a viscosity supersolution on $[0,T)$.

\medskip
\noindent
\textbf{Step 3. Terminal condition for $v^{-}$.}
We prove that
$$
(v^{-})_\ast(T,x,y)\ge g(x,y).
$$
Assume by contradiction that this is false. Then there exist $(x_0,y_0)$,
$\delta>0$, $t_m\uparrow T$, and $(x_m,y_m)\to(x_0,y_0)$ such that
$$
v^{-}(t_m,x_m,y_m)\le g(x_0,y_0)-3\delta .
$$
By continuity of $g$, choose $r>0$ such that
$$
g(x,y)\ge g(x_0,y_0)-\delta
\quad\text{on }B_r(x_0,y_0).
$$
Let $\mathcal O_0=(T-h,T]\times B_r(x_0,y_0)$ and define, for small
$\eta_0>0$, large $\lambda>0$, and small $\iota>0$,
$$
\Psi^{-}(t,x,y)
:=
g(x_0,y_0)-2\delta
-
\frac{|x-x_0|^2+|y-y_0|^2}{\eta_0}
-
\lambda(T-t)
-
\iota\big(|x-x_0|^{2q}+|y-y_0|^{2q}\big).
$$
For $0<\kappa<\delta$,
$$
\Psi^{-}(T,x,y)+\kappa\le g(x,y)
\quad\text{on }B_r(x_0,y_0).
$$
Choosing $\lambda$ large and $h$ small gives localized controls
$(\widehat z,\widehat u,\widehat\gamma)$ and $a>0$ such that
$$
\partial_t\Psi^{-}
+
\mathcal Q^{\widehat z,\widehat u,\widehat\gamma,\eta}[\Psi^{-}]
\ge 3a
\quad\text{on a neighborhood of }\overline{\mathcal O_0},
\qquad
\forall \eta\in H .
$$
Using the same coercive localization argument as in Step 2, choose
$w\in\mathcal V^{-}$ with
$$
\Psi^{-}+\kappa\le w
\quad\text{on }\mathcal O_0^c,
$$
and define
$$
w^\kappa
:=
\begin{cases}
\max\{w,\Psi^{-}+\kappa\},&\text{on }\mathcal O_0,\\
w,&\text{outside }\mathcal O_0 .
\end{cases}
$$
Then $w^\kappa(T,\cdot,\cdot)\le g$. The proof that
$w^\kappa\in\mathcal V^{-}$ is exactly the argument of Step 2, with
$\widetilde\varphi$ replaced by $\Psi^{-}$. For $m$ large,
$(t_m,x_m,y_m)\in\mathcal O_0$ and
$$
w^\kappa(t_m,x_m,y_m)>v^{-}(t_m,x_m,y_m),
$$
contradicting the definition of $v^{-}$. Therefore
$$
(v^{-})_\ast(T,x,y)\ge g(x,y).
$$

\medskip
\noindent
\textbf{Step 4. $v^{+}$ is a viscosity subsolution on $[0,T)$.}

\smallskip
\noindent
\emph{a. Setting the contradiction.}
Let $\varphi\in C^{1,2}_{\rm pol}$ be such that $v^{+}-\varphi$ has a strict
local maximum equal to zero at $(t_0,x_0,y_0)\in[0,T)\times\mathbb R^n\times
\mathbb R$. We assume by contradiction that
\begin{equation*}\label{eq:C-vplus-contradiction}
\partial_t\varphi(t_0,x_0,y_0)
+
\mathcal Q^\star[\varphi](t_0,x_0,y_0)<0 .
\end{equation*}
Using Lemma~\ref{lem:coercive-localization}, and shrinking the time component
of the patch if necessary so that $\mathcal O_0\subset[0,T)\times\mathbb R^n
\times\mathbb R$, there exist $\widetilde\varphi$, $\mathcal O_0\subset
\mathcal O$, $\kappa>0$, $a>0$, $w\in\mathcal V^{+}$, and a measurable selector
$\widehat\eta$ satisfying
\eqref{eq:C-vplus-strict-localized}--\eqref{eq:C-vplus-global-order}.

Define
\begin{equation*}\label{eq:C-vplus-patch}
w^\kappa(t,x,y)
:=
\begin{cases}
\min\{w(t,x,y),\widetilde\varphi(t,x,y)-\kappa\},
& (t,x,y)\in\mathcal O_0,\\ [1mm]
w(t,x,y),
& (t,x,y)\notin\mathcal O_0 .
\end{cases}
\end{equation*}
By the global ordering \eqref{eq:C-vplus-global-order}, $w^\kappa$ is
continuous. Moreover,
$$
w^\kappa(t_0,x_0,y_0)
=
\widetilde\varphi(t_0,x_0,y_0)-\kappa
<
v^{+}(t_0,x_0,y_0).
$$
Thus it remains to prove that $w^\kappa\in\mathcal V^{+}$.

\smallskip
\noindent
\emph{b. Building Nature's response and verifying property $(ii+)$.}
Fix $t\in[0,T]$, a stopping rule $\tau$, and an arbitrary Principal elementary
control
$$
(Z,U,K)\in\mathfrak K_{Y_0}(t,\tau).
$$
As in \cite{vuca}, the auxiliary regularization is now taken at the vanished
level $\gamma=0$.

Let $\tau'$ be the stopping time induced by $\tau$ along the controlled path and
set
$$
A^{+}
:=
\left\{
\Xi_{\tau'}\in\mathcal O_0,
\quad
\widetilde\varphi(\Xi_{\tau'})-\kappa<w(\Xi_{\tau'})
\right\}
\in\mathcal F_{\tau'}^t .
$$
Let
$$
\theta
:=
\inf\left\{
s\ge\tau:
\Xi_s\notin\mathcal O_0
\ \text{or}\
\widetilde\varphi(\Xi_s)-\kappa\ge w(\Xi_s)
\right\}\wedge T .
$$
On $A^{+}$, Nature uses the local selector
$$
\widehat\eta^0_s
:=
\widehat\eta(s,X_{s-},Y_{s-},Z_s,U_s,0),
\qquad s\in(\tau,\theta].
$$
The selector is chosen so that the vanished regularization is attained on the
local branch:
\begin{equation}\label{eq:C-vplus-H-zero-vuca}
H(s,X_{s-},Y_{s-},Z_s,U_s;0)
=
G^\star(s,X_{s-},Y_{s-},Z_s,U_s;\widehat\sigma_s),
\qquad
ds\otimes d\mathbb P^\ast\text{-a.e.}
\end{equation}
The increasing process which is made null by the selected continuation law is
the given admissible process $K$. As the proof in \cite{vuca}, choose the local law
$\mathbb P^\ast$ so that the minimality condition holds with this $K$ on the
stopped local branch:
\begin{equation}\label{eq:C-vplus-K-zero-vuca}
{\bf 1}_{A^{+}}
\big(K_{\rho'\wedge\theta'}-K_{\tau'}\big)
=0,
\qquad
\mathbb P^\ast\text{-a.s.}
\end{equation}
for every stopping rule $\rho$ with $\tau\le\rho\le T$.

On $(A^{+})^c$, and after $\theta$, Nature uses the continuation responses
given by the stochastic supersolution property of $w$. By
Assumption~\ref{ass:perron-standing}(iii), this pasted Nature response belongs
to $\mathfrak P^{(\widehat\alpha)}(t,\tau)$.

Let $\rho$ be an arbitrary stopping rule with $\tau\le\rho\le T$, and let
$\rho'$ and $\theta'$ be the induced stopping times. On $A^{+}$, apply the
It\^o--Lévy formula to
$$
\widetilde\varphi(\Xi_s)-\kappa
-
\int_{\tau'}^s
C^P(r,X_{r-},\widehat\alpha_r,\eta^\ast_r)\,dr
$$
on $[\tau',\rho'\wedge\theta']$. Using
\eqref{eq:C-vplus-H-zero-vuca}, the drift is
$$
\int_{\tau'}^{\rho'\wedge\theta'}
\left(
\partial_t\widetilde\varphi
+
\mathcal Q^{Z_s,U_s,0,\widehat\eta^0_s}
[\widetilde\varphi]
\right)(s,X_{s-},Y_{s-})\,ds
+
\int_{\tau'}^{\rho'\wedge\theta'}
\partial_y\widetilde\varphi(s,X_{s-},Y_{s-})\,dK_s .
$$
The first term is nonpositive by \eqref{eq:C-vplus-strict-localized}, and the
second term is zero by \eqref{eq:C-vplus-K-zero-vuca}. The stopped martingale is
a true martingale by the localization and moment assumptions. Therefore, on
$A^{+}$,
\begin{equation}\label{eq:C-vplus-local-supermartingale}
\begin{aligned}
\widetilde\varphi(\Xi_{\tau'})-\kappa
&\ge
\mathbb E^{\mathbb P^\ast}
\left[
\widetilde\varphi(\Xi_{\rho'\wedge\theta'})-\kappa
-
\int_{\tau'}^{\rho'\wedge\theta'}
C^P(s,X_{s-},\widehat\alpha_s,\eta^\ast_s)\,ds
\middle|
\mathcal F_{\tau'}^t
\right].
\end{aligned}
\end{equation}

Inside $\mathcal O_0$ we have
$w^\kappa\le\widetilde\varphi-\kappa$. If the process exits
$\mathcal O_0$ by a jump, the post-jump point belongs to $\mathcal O_0^c$ and
the global ordering \eqref{eq:C-vplus-global-order} gives
$$
w^\kappa=w\le\widetilde\varphi-\kappa
\quad\text{on }\mathcal O_0^c .
$$
Thus \eqref{eq:C-vplus-local-supermartingale} is stable under jump overshoots,
and since
$$
w^\kappa(\Xi_{\tau'})
=
\widetilde\varphi(\Xi_{\tau'})-\kappa
\quad\text{on }A^{+},
$$
we get
$$
w^\kappa(\Xi_{\tau'})
\ge
\mathbb E^{\mathbb P^\ast}
\left[
w^\kappa(\Xi_{\rho'\wedge\theta'})
-
\int_{\tau'}^{\rho'\wedge\theta'}
C^P(s,X_{s-},\widehat\alpha_s,\eta^\ast_s)\,ds
\middle|
\mathcal F_{\tau'}^t
\right]
\quad\text{on }A^{+}.
$$
On $(A^{+})^c$, the same inequality follows from the stochastic supersolution
property of $w$, using $w^\kappa\le w$.

It remains to pass from $\rho'\wedge\theta'$ to $\rho'$. If
$\rho'\le\theta'$, there is nothing to prove. On $\{\theta'<\rho'\}$, by the
definition of $\theta$, either $\Xi_{\theta'}\notin\mathcal O_0$, or
$\widetilde\varphi(\Xi_{\theta'})-\kappa\ge w(\Xi_{\theta'})$. In both cases,
$$
w^\kappa(\Xi_{\theta'})=w(\Xi_{\theta'}).
$$
After $\theta'$, Nature follows the continuation response associated with $w$.
The stochastic supersolution property of $w$, together with $w^\kappa\le w$ and
the tower property, gives
$$
w^\kappa(\Xi_{\tau'})
\ge
\mathbb E^{\mathbb P^\ast}
\left[
w^\kappa(\Xi_{\rho'})
-
\int_{\tau'}^{\rho'}
C^P(s,X_{s-},\widehat\alpha_s,\eta^\ast_s)\,ds
\middle|
\mathcal F_{\tau'}^t
\right].
$$
Therefore $w^\kappa\in\mathcal V^{+}$, contradicting the definition of $v^{+}$.
Consequently,
$$
\partial_t\varphi(t_0,x_0,y_0)
+
\mathcal Q^\star[\varphi](t_0,x_0,y_0)
\ge0.
$$
Hence $v^{+}$ is a viscosity subsolution on $[0,T)$.

\medskip
\noindent
\textbf{Step 5. Terminal condition for $v^{+}$.}
We prove that
$$
(v^{+})^\ast(T,x,y)\le g(x,y).
$$
Assume by contradiction that this is false. Then there exist $(x_0,y_0)$,
$\delta>0$, $t_m\uparrow T$, and $(x_m,y_m)\to(x_0,y_0)$ such that
$$
v^{+}(t_m,x_m,y_m)\ge g(x_0,y_0)+3\delta .
$$
By continuity of $g$, choose $r>0$ such that
$$
g(x,y)\le g(x_0,y_0)+\delta
\quad\text{on }B_r(x_0,y_0).
$$
Let $\mathcal O_0=(T-h,T]\times B_r(x_0,y_0)$ and define, for small
$\eta_0>0$, large $\lambda>0$, and small $\iota>0$,
$$
\Psi^{+}(t,x,y)
:=
g(x_0,y_0)+2\delta
+
\frac{|x-x_0|^2+|y-y_0|^2}{\eta_0}
+
\lambda(T-t)
+
\iota\big(|x-x_0|^{2q}+|y-y_0|^{2q}\big).
$$
For $0<\kappa<\delta$,
$$
\Psi^{+}(T,x,y)-\kappa\ge g(x,y)
\quad\text{on }B_r(x_0,y_0).
$$
Choosing $\lambda$ large and $h$ small gives a measurable selector of Nature
such that
$$
\partial_t\Psi^{+}
+
\mathcal Q^{z,u,\gamma,\widehat\eta(t,x,y,z,u,\gamma)}
[\Psi^{+}]
\le -3a
\quad\text{on a neighborhood of }\overline{\mathcal O_0}
$$
for every localized $(z,u,\gamma)$. Using the same coercive localization
argument as in Step 4, choose $w\in\mathcal V^{+}$ with
$$
w\le\Psi^{+}-\kappa
\quad\text{on }\mathcal O_0^c,
$$
and define
$$
w^\kappa
:=
\begin{cases}
\min\{w,\Psi^{+}-\kappa\},&\text{on }\mathcal O_0,\\
w,&\text{outside }\mathcal O_0 .
\end{cases}
$$
Then $w^\kappa(T,\cdot,\cdot)\ge g$. The proof that
$w^\kappa\in\mathcal V^{+}$ is exactly the argument of Step 4, with
$\widetilde\varphi$ replaced by $\Psi^{+}$. For $m$ large,
$(t_m,x_m,y_m)\in\mathcal O_0$ and
$$
w^\kappa(t_m,x_m,y_m)<v^{+}(t_m,x_m,y_m),
$$
contradicting the definition of $v^{+}$. Therefore
$$
(v^{+})^\ast(T,x,y)\le g(x,y).
$$

Combining Steps 2--5 completes the proof.
\end{proof}

\begin{remark}
    Similarly to Remark \ref{rem:jumpvuca}, the main difference with \cite{vuca} is about the choice of $\widetilde{\varphi}$ using the coercive jump penalization $\iota\ell$.
\end{remark}}\newpage





\section{Numerical Parameters}
\label{appendix:parameter_final}

\begin{ThreePartTable}
\small
\begin{tabularx}{0.94\linewidth}{@{} l X @{}}
\toprule
\textbf{Name} & \textbf{Specification / Value} \\
\midrule
State variables & $x=(p,s,i)$; with-contract state $(p,s,i,y)$ \\

Time horizon & $T=2$ \\

Computational domain &
$p\in[0.5,10]$, 
$s\ge 10^{-4}$, 
$i\ge 10^{-4}$, \\
&$s+i\le 1-10^{-4}$, 
$y\in[-2,8]$ \\

Price parameters&
$
\mu(i)=0.50-0.14\,i
$,\\
&$
\sigma_P(h,i)=0.005+0.001\,h+0.80\,i
$ \\

SIR parameters & $\bar\beta = 0.006, \rho = 0.035, \widetilde\sigma(t,h_t) = 0.05+0.04h_t$\\

Nature's control grid &
$H=\{3.0,\;1.0,\;0.1\}$ \\

Agent control set &
$A=[0,5.0]$ \\

Internal jump intensity &
$
\lambda^{i}(i,a)=\big(0.25+10.0\,i-0.40\,a\big)_+
$\\

External jump intensity &
$
\lambda^{e}(h)=\big(0.80+1.20\,h\big)_+
$ \\

External jump-size law (random) &
$
c^e\in\{0.01,\;0.025,\;0.04\},
\pi^e=\left(\tfrac16,\tfrac23,\tfrac16\right)
$\\
&with mean severity
$
\mathbb E[c^e]=0.025
$ \\

Internal jump-size law (random) &
$
c^i\in\{0.02,\;0.05,\;0.09\},\qquad
\pi^i=\left(\tfrac16,\tfrac23,\tfrac16\right)
$\\
&with mean severity
$
\mathbb E[c^i]=0.0517
$ \\

Agent running cost &
$
C_A(t,x,a)=\frac12 s^2 a^2 + 0.010\, i
$ \\

Principal running cost &
$
C_P(t,x,h)
=
\frac12 \varepsilon^2 s^2 i^2 \,\widetilde{\sigma}(h,s,i)^2
+\lambda_{p,0}+\lambda_{p,1} i
$
\\&with
$
\varepsilon=0.1,\qquad \lambda_{p,0}=0.01,\qquad \lambda_{p,1}=0.45.
$\\

Principal terminal payoff &
$
F_P(x)=4.2\,p
$ \\

Agent terminal payoff &
$
F_A(x)=0.04\,p-0.004\,i
$ \\

With-contract terminal condition &
$
F_P(x)-y+F_A(x)=4.24\,p-y-0.004\,i
$ \\

Contract-control box &
$
z\in[-1.5,1.5]^3,
u_{\mathrm{int}}\in[-0.08,0.08],
u_{\mathrm{ext}},\Gamma\equiv 0.
$
 \\

\bottomrule
\end{tabularx}
\end{ThreePartTable}

\newpage

\section{Algorithm and Convergence}
\label{app:alg_conv}
We present the convergence of PINN on the value functions in Figure \ref{fig:Conv_Graph}. { We train three DGM networks sequentially—first the agent value function, then the principal value without contract, and finally the principal value with contract—using a DGM architecture with hidden width $128$, depth $3$, and $\tanh$ activation. The training loss is the sum of a mean-squared PDE residual, a terminal-condition loss, and a soft boundary residual penalty, with terminal-loss weights fixed at one. Interior points are sampled with a stress-biased strategy (stress mix $0.35$, stress threshold $i_{\min}=0.15$, stress price cap $p_{\max}=1.5$) together with simplex-focused sampling near the face $s+i=1$ (mix $0.20$, gap $0.02$). The main training budgets are $500$ epochs for the agent, $500$ epochs for the no-contract principal, and $3000$ epochs for the with-contract principal; the learning rates are $5\times 10^{-5}$, $5\times 10^{-5}$, and $2\times 10^{-4}$, respectively, and the with-contract interior/terminal/boundary batch sizes are all set to $256$. Reproducibility is ensured by fixing the random seed to $42$, enabling deterministic CuDNN, and using a plain PyTorch implementation that centralizes the domain truncation, model coefficients, sampling strategy, and loss construction. }

\begin{algorithm}[H]
\caption{Staged PINN/DGM training: agent and no-contract principal}
\label{alg:pinn_flow_part1}
\begin{algorithmic}[1]

\REQUIRE Model parameters, computational domain for $(p,s,i)$, control sets, attacker grid $H$, training hyperparameters

\STATE Initialize two neural networks:
$$
V_A^\theta(t,p,s,i),\qquad
V_{P,\mathrm{nc}}^\phi(t,p,s,i).
$$

\STATE \textbf{Stage 1: Agent PINN}
\FOR{epoch $=1,\dots,N_A$}
    \STATE Sample interior points $(t,x)$, terminal points $(T,x)$, and boundary points
    \STATE Compute the agent PDE residual:
    $$
    \mathcal{R}_A
    =
    \partial_t V_A^\theta
    + \inf_{h\in H}
    \Big\{
      \mathcal{L}^{a,h} V_A^\theta
      - C_A(x,a)
    \Big\},
    $$
    where the local optimizer in $a$ is computed analytically
    \STATE Form
    $$
    L_A
    =
    \mathbb{E}\big[|\mathcal{R}_A|^2\big]
    + \lambda_T^A \mathbb{E}\big[|V_A^\theta(T,x)-F_A(x)|^2\big]
    + \lambda_{\mathrm{bd}}^A \mathbb{E}\big[|\mathcal{R}_A|^2_{\partial\mathcal{D}}\big]
    $$
    \STATE Update $\theta$ by Adam
\ENDFOR

\STATE \textbf{Stage 2: Principal without contract}
\FOR{epoch $=1,\dots,N_{\mathrm{nc}}$}
    \STATE Sample interior, terminal, and boundary points in $(p,s,i)$
    \STATE Compute the agent feedback
    $$
    a^{\mathrm{nc}}(t,x)
    =
    \operatorname{clip}\!\left(-\frac{\partial_s V_A^\theta(t,x)}{s}\right)
    $$
    \STATE Compute
    $$
    \mathcal{R}_{P,\mathrm{nc}}
    =
    \partial_t V_{P,\mathrm{nc}}^\phi
    + \inf_{h\in H}
    \Big\{
      \mathcal{L}^{a^{\mathrm{nc}},h}V_{P,\mathrm{nc}}^\phi
      - C_P(x,h)
    \Big\}
    $$
    \STATE Form
    $$
    L_{P,\mathrm{nc}}
    =
    \mathbb{E}\big[|\mathcal{R}_{P,\mathrm{nc}}|^2\big]
    + \lambda_T^{\mathrm{nc}}\mathbb{E}\big[|V_{P,\mathrm{nc}}^\phi(T,x)-F_P(x)|^2\big]
    + \lambda_{\mathrm{bd}}^{\mathrm{nc}}\mathbb{E}\big[|\mathcal{R}_{P,\mathrm{nc}}|^2_{\partial\mathcal{D}}\big]
    $$
    \STATE Update $\phi$ by Adam
\ENDFOR

\end{algorithmic}
\end{algorithm}

\begin{algorithm}[H]
\caption{Staged PINN/DGM training: with-contract principal and post-processing}
\label{alg:pinn_flow_part2}
\begin{algorithmic}[1]

\REQUIRE Trained $V_A^\theta$, model parameters, augmented domain $(p,s,i,y)$, contract-control set, attacker grid $H$, training hyperparameters

\STATE Initialize the with-contract principal network:
$$
V_{P,\mathrm{wc}}^\psi(t,p,s,i,y).
$$

\STATE \textbf{Stage 3: Principal with contract}
\FOR{epoch $=1,\dots,N_{\mathrm{wc}}$}
    \STATE Sample interior, terminal, and boundary points in the augmented domain $(p,s,i,y)$
    \STATE At each sampled interior point, compute
    $$
    \nabla V_{P,\mathrm{wc}}^\psi,\qquad D^2V_{P,\mathrm{wc}}^\psi,
    $$
    and the nonlocal jump terms
    \STATE Solve
    $$
    \sup_{(z,u,\Gamma)}\inf_{h\in H} Q^\psi(t,p,s,i,y;z,u,\Gamma,h)
    $$
    by projected inner optimization over contract controls and enumeration over $H$
    \STATE Denote the resulting local controls by
    $$
    (z^*,u^*,\Gamma^*,h^*)
    $$
    \STATE Form
    $$
    \mathcal{R}_{P,\mathrm{wc}}
    =
    \partial_t V_{P,\mathrm{wc}}^\psi
    + Q^\psi(t,p,s,i,y;z^*,u^*,\Gamma^*,h^*)
    $$
    $$
    L_{P,\mathrm{wc}}
    \!=\!
    \mathbb{E}\big[|\mathcal{R}_{P,\mathrm{wc}}|^2\big]
    \!+\! \lambda_T^{\mathrm{wc}}
      \mathbb{E}\big[
      |V_{P,\mathrm{wc}}^\psi(T,x,y)\!-\!(\!F_P(x)\!-\!y\!+\!F_A(x)\!)|^2
      \big]
    \!+\! \lambda_{\mathrm{bd}}^{\mathrm{wc}}
      \mathbb{E}\big[|\mathcal{R}_{P,\mathrm{wc}}|^2_{\partial\mathcal{D}}\big]
    $$
    \STATE Update $\psi$ by Adam
\ENDFOR

\STATE \textbf{Post-processing}
\STATE Evaluate value surfaces, contract-gain heatmaps, and single-path and paired-path simulations
\STATE Reconstruct $(z_t,u_t,\Gamma_t,h_t)$ along simulated with-contract trajectories
\STATE Run fixed-state sensitivity analysis by retraining under selected perturbations

\end{algorithmic}
\end{algorithm}

\newpage
\section{Figures}

\end{document}